\documentclass[letterpaper, 10 pt, conference]{ieeeconf}  

\newcommand{\ls}[2]{#2}\renewcommand{\baselinestretch}{1}

\IEEEoverridecommandlockouts                              
\usepackage[english]{babel}
\usepackage[T1]{fontenc}
\usepackage{graphicx} 
\graphicspath{{figures/}} 
\usepackage{amsmath}
\usepackage{amsfonts}
\usepackage{mathtools}
\usepackage{tensor}

\usepackage{amssymb}
\usepackage{arydshln}
\usepackage{bm}
\usepackage{cite}
\usepackage{tikz,datatool}
\usepgflibrary{fpu}
\usetikzlibrary{positioning,calc,arrows,shapes,shadows,matrix,backgrounds}
\usepackage{xcolor}
\usepackage{hyperref}

\usetikzlibrary{positioning,calc,arrows,shapes,shadows}
\tikzset{
auto,
sys/.style 2 args={
rectangle,
draw,
drop shadow,
fill=white,
minimum height=#2,
minimum width=#1,
inner sep=\dn},
sum/.style={circle,draw,draw=black,inner sep=0mm,minimum size=2mm},
jun/.style={circle,draw,draw=black,inner sep=0mm,minimum size=0mm},
>={latex},
every path/.style={rounded corners},
}

\newcommand{\tio}[4]{\coordinate (#1) at ($(#2.south #3)!#4!(#2.north #3)$)}

\def\dn{1ex}
\def\dl{3*\dn}

\tikzstyle{sy0}=[sys={0*\dn}{0*\dn}]
\tikzstyle{sy1}=[sys={12*\dn}{8*\dn}]
\tikzstyle{sy2}=[sys={8*\dn}{6*\dn}]
\tikzstyle{sy3}=[sys={5*\dn}{5*\dn}]

\newcommand{\arc}{\arraycolsep0.3ex}
\newcommand{\bul}{\bullet}
\newcommand{\im}{{\rm im}}

\newcommand{\Af}{A_{\rm f}}
\newcommand{\Bf}{B_{\rm f}}

\newcommand{\Vf}{V_{\rm f}}
\newcommand{\Va}{ {{\hat V_2}} }
\newcommand{\Vb}{ {{\hat V_3}} }
\newcommand{\Cfla}{\bm C_{\rm\bf f}}
\newcommand{\Dfla}{\bm D_{\rm\bf f}}
\newcommand{\Mla}{\bm M}
\newcommand{\Nla}{\bm N}
\newcommand{\Ttla}{\bm {\tilde T}}
\newcommand{\Tla}{\bm T}

\newcommand{\Yf}{Y_{\rm f}}

\newcommand{\Tf}{T_{\rm f}}

\newcommand{\At}{\tilde A}
\newcommand{\Bt}{{\tilde B}}
\newcommand{\Ct}{\tilde C}
\newcommand{\Ctz}{\tilde C_z}
\newcommand{\Dt}{\tilde D}
\newcommand{\Et}{{\tilde D}_z}
\newcommand{\Yt}{\tilde Y}

\newcommand{\Yh}{\hat Y}
\newcommand{\Ah}{\h A}
\newcommand{\Bh}{\h B}
\newcommand{\Ch}{{\h C}}
\newcommand{\Dh}{{\h D}}

\newcommand{\Uh}{\h U}
\newcommand{\Vh}{\h V}

\newcommand{\Chla}{\bm \Ch_z}
\newcommand{\Dhla}{\bm \Dh_{zw}}
\newcommand{\Ehla}{\bm \Dh_z}

\DeclareMathAlphabet\mathbfcal{OMS}{cmsy}{b}{n}

\newcommand{\Acl}{{\cal A}}
\newcommand{\Bcl}{{\cal B}}

\newcommand{\Ccla}{\mathbfcal{C}}
\newcommand{\Dcla}{\mathbfcal{D}}

\newcommand{\Xcl}{{{\cal X} }}
\newcommand{\Xclf}{{\cal X_{\rm f}}}
\newcommand{\Xcls}{{\cal X_{\rm s}}}
\newcommand{\Xclfs}{{\cal X_{\rm fs}}}
\newcommand{\Xclsf}{{\cal X_{\rm sf}}}

\newcommand{\Al}{A}
\newcommand{\Bl}{B}
\newcommand{\Cl}{C}
\newcommand{\Dl}{D}
\newcommand{\Gl}{G}
\newcommand{\all}{\al}
\newcommand{\bel}{\be}
\newcommand{\nl}{n}
\newcommand{\Sl}{S}
\newcommand{\Tl}{T}

\newcommand{\Ar}{\h A}
\newcommand{\Br}{\h B}
\newcommand{\Cr}{\h C}
\newcommand{\Dr}{\h D}
\newcommand{\Gr}{\hat G}
\newcommand{\alr}{\hat\al}

\newcommand{\nr}{\hat n}

\newcommand{\Id}{I_d}
\newcommand{\ot}{\otimes}
\newcommand{\ots}{\!\otimes\!}

\newcommand{\arr}[2]{\begin{array}{#1}#2\end{array}}

\newcommand{\te}[1]{\text{\ \ #1\ \ }}

\newcommand{\z}{{\rm z}}

\renewcommand{\r}[1]{(\ref{#1})}
\newcommand{\hl}{\\\hline}

\newcommand{\enu}[1]{\begin{enumerate}#1\end{enumerate}}
\newcommand{\mun}[1]{\vspace*{-1ex}\begin{multline*}#1\end{multline*}}
\newcommand{\mul}[1]{\vspace*{-1ex}\begin{multline}#1\end{multline}}
\newcommand{\mas}[2]{\left[\begin{array}{#1}#2\end{array}\right]}
\newcommand{\smat}[2]{{\tiny \left(\begin{array}{#1}#2\end{array}\right)}}
\renewcommand{\c}[1]{{\cal #1}}
\newcommand{\eql}[2]{\begin{equation}\label{#1}#2\end{equation}}
\newcommand{\T}{\top}

\newcommand{\Aalg}{A_a}
\newcommand{\Balg}{B_a}
\newcommand{\Calg}{C_a}

\newcommand{\Aalp}{{A_c}}
\newcommand{\Balp}{{B_c}}
\newcommand{\Calp}{{C_c}}
\newcommand{\Dalp}{{D_c}}

\newcommand{\Aco}{{A_c}}
\newcommand{\Bco}{{B_c}}
\newcommand{\Cco}{{C_c}}
\newcommand{\Dco}{{D_c}}

\newcommand{\Ach}{{\hat A_c}}
\newcommand{\Bch}{{\hat B_c}}
\newcommand{\Cch}{{\hat C_c}}
\newcommand{\Dch}{{\hat D_c}}

\newcommand{\As}{A}
\newcommand{\Bsw}{B_w}
\newcommand{\Bs}{B}
\newcommand{\Csz}{C_z}
\newcommand{\Dsz}{D_z}
\newcommand{\Cs}{C}

\newcommand{\R}{\mathbb{R}}         
\newcommand{\N}{\mathbb{N}}         
\renewcommand{\S}{\mathbb{S}}          
\newcommand{\F}{\mathbb{F}}

\newcommand{\I}{\mathbb{I}}
\newcommand{\C}{\mathbb{C}}
\newcommand{\D}{\mathbb{D}}

\newcommand{\cle}{\preccurlyeq}
\newcommand{\cl}{\prec}
\newcommand{\cg}{\succ}

\newcommand{\diag}{\mathrm{diag}}
\newcommand{\col}{\mathrm{col}}

\newcommand{\ga}{\gamma}
\newcommand{\eps}{\varepsilon}
\newcommand{\la}{\lambda}
\newcommand{\La}{\Lambda}

\newcommand{\al}{\alpha}
\newcommand{\be}{\beta}

\newcommand{\mat}[2]{\left(\begin{array}{@{}#1@{}}#2\end{array}\right)} 

\newenvironment{red_test}{\color{red}}{} 

\newenvironment{blue_test}{\color{blue}}{} 

\renewcommand{\t}{\widetilde}
\newcommand{\h}{\hat}

\DeclareMathSymbol{\shortminus}{\mathbin}{AMSa}{"39}

\newcommand{\mstrut}[1]{\rule{0pt}{#1}}
\let\oldhline=\hline 
\renewcommand{\hline}{\oldhline\mstrut{2.25ex}}
\let\oldhdashline=\hdashline 
\renewcommand{\hdashline}{\oldhdashline\mstrut{2.25ex}}

\newtheorem{theorem}{Theorem}
\newtheorem{lemma}[theorem]{Lemma}

\newtheorem{definition}[theorem]{Definition}

\newcounter{frstequation}
\newcounter{scndequation}
\newcounter{trdequation}

\makeatletter
\newcommand{\scndlabel}[1]{%
	\protected@write \@auxout {}{\string \newlabel{#1}{{\@currentlabel}{\thepage}{}{equation.\thesection.\the\value{parentequation}\alph{frstequation}}{}} }%
	\hypertarget{equation.\thesection.\the\value{parentequation}\alph{frstequation}}{}
}
\makeatother
\makeatletter
\newcommand{\dlabel}[2]{ 
	\setcounter{frstequation}{\value{equation}}%
	\setcounter{scndequation}{\value{equation}}%
	\addtocounter{scndequation}{1}%
	\edef\theequation{\theparentequation\alph{frstequation},\alph{scndequation}}
	\def\@currentlabel{\theparentequation\alph{frstequation}}\scndlabel{#1}
	\def\@currentlabel{\theparentequation\alph{scndequation}}\label{#2}
	\addtocounter{equation}{1}
}
\makeatother

\makeatletter
\newcommand{\tlabel}[3]{ 
	\setcounter{frstequation}{\value{equation}}
	\setcounter{scndequation}{\value{equation}}%
	\setcounter{trdequation}{\value{equation}}%
	\addtocounter{scndequation}{1}%
	\addtocounter{trdequation}{2}%
	\edef\theequation{\theparentequation\alph{frstequation},\alph{scndequation},\alph{trdequation}}
	\def\@currentlabel{\theparentequation\alph{frstequation}}\scndlabel{#1}
	\def\@currentlabel{\theparentequation\alph{scndequation}}\scndlabel{#2}
	\def\@currentlabel{\theparentequation\alph{trdequation}}\label{#3}
	\addtocounter{equation}{2}
}
\makeatother

\def\BibTeX{{\rm B\kern-.05em{\sc i\kern-.025em b}\kern-.08em
		T\kern-.1667em\lower.7ex\hbox{E}\kern-.125emX}}
\newcommand{\epro}{\mbox{}\hfill\mbox{\rule{2mm}{2mm}} }

\title{Optimization Algorithm Synthesis based on Integral Quadratic Constraints: A Tutorial}

\author{Carsten W. Scherer,  Christian Ebenbauer and Tobias Holicki
\thanks{The first and third author are funded by Deutsche Forschungsgemeinschaft (DFG, German Research Foundation) under Germany's Excellence Strategy - EXC 2075 - 390740016. They acknowledge the support by the Stuttgart Center for Simulation Science (SimTech).}
\thanks{Carsten W. Scherer and Tobias Holicki are with the Department of Mathematics, University of Stuttgart, Germany (email: \{carsten.scherer,tobias.holicki\}@imng.uni-stuttgart.de).
Christian Ebenbauer is  the Chair of Intelligent Control Systems, RWTH Aachen University, Germany (email: christian.ebenbauer@ic.rwth-aachen.de).
}
}

\begin{document}

\maketitle

\begin{abstract}
We expose in a tutorial fashion the mechanisms which underlie the synthesis of optimization algorithms based on dynamic integral quadratic constraints.
We reveal how these tools from robust control allow to
design accelerated gradient descent algorithms with optimal guaranteed convergence rates by solving small-sized convex semi-definite programs. It is shown that this extends to the design of extremum controllers, with the goal to
regulate the output of a general linear closed-loop system to the minimum of an objective function.

Numerical experiments illustrate that we can not only recover gradient decent and the triple momentum variant of Nesterov's accelerated first order algorithm, but also automatically synthesize optimal algorithms even if the gradient information
is passed through non-trivial dynamics, such as time-delays.
\end{abstract}




\section{Introduction}

Accelerated gradient algorithms \cite{Nes18} have a wide range of applications in the current era of
machine learning and online optimization-based control. From the perspective of control theory,
such algorithms can be viewed as a linear time-invariant discrete-time (LTI) system in feedback with the
gradient of the to-be-minimized function as a nonlinearity
\cite{Pol87,WanEli11,DueEbe12,LesRec16}. This provides an immediate link to
absolute stability theory and offers the possibility to apply advanced
tools from robust control for the {\em automated analysis} of accelerated gradient algorithms \cite{LesRec16}.
By tuning the algorithm parameters based on these tools, the convergence rate of Nesterov’s algorithm \cite{Nes83} has been improved to get the so-called triple momentum algorithm \cite{ScoFre18}.

The {\em automated synthesis} of optimization algorithms by convex optimization
is a much more challenging task. This falls into the area of robust feedback controller design \cite{ZhoDoy96,SchWei11}.
Recent work \cite{LesSei20,MicSch21,GraEbe22} has addressed the synthesis problem from this perspective,
but based on heuristic methods without optimality guarantees.
An alternative approach to non-convex algorithm design by interpolation techniques can be found in \cite{TayHen16,TayDro22}.

The purpose of this paper is to develop, in a tutorial fashion, the whole pipeline
of analysis techniques that open the avenue for a convex
solution to the automated algorithm synthesis problem by solving a
moderate-sized convex semi-definite program. Another feature of the presented approach
is its flexibility. It offers a convex solution to the so-called extremum control problem, with the goal to regulate the output of a dynamical system to
the minimum of some convex cost function. These main results are based on \cite{SchEbe21,HolSch21c}.
However, we also present an innovation over \cite{HolSch21c} which renders synthesis possible
for LTI systems without any restrictions on their poles or zeros.

The paper is structured as follows. In Sec. \ref{SecII}, we show how the
algorithm analysis and synthesis problems translate into one of robustness analysis and synthesis.
Sec. \ref{SecIII} recaps robustness analysis with static integral quadratic
constraints (IQCs). Dynamic IQCs are introduced in Sec. \ref{SecIV}, while the corresponding robust stability test is given in Sec. \ref{SecV}.
The design of algorithms is presented in Sec.~\ref{Ssyn} and numerical illustrations
are found in Sec. \ref{SNum}. Concluding remarks are given in Sec. \ref{Scon}.
All proofs and some explanatory connections to classical passivity-based stability tests
are found in \ls{an extended version of the paper on arXiv \cite{SchEbe23a}.}{the appendix.}

Next to standard notations, for matrices $A$, $B$ we express by $A\leq B$ that $B-A$ is nonnegative entrywise,
while $A\cl B$ means that $A$ and $B$ are symmetric and $B-A$ is positive definite.
For a tuple of matrices $A=(A_1,\ldots,A_k)$, we use
$$\renewcommand{\arraystretch}{.9}
\diag(A)=\mat{ccc}{A_1&\cdots&0\\[-1ex]\vdots&\ddots&\vdots\\0&\cdots&A_k}\te{and}
\col(A)=\mat{c}{A_1\\[-1ex]\vdots\\A_k}
$$
if the dimensions are compatible.
For the real polynomial $\al(\z)=\al_0+\cdots+\al_{n-1}\z^{n-1}+\z^n$ of degree $n$,
we denote by $C_\al\in\R^{n\times n}$ the standard companion matrix
with the last row $(-\al_0,\ldots,-\al_{n-1})$, and $e_n\in\R^n$ is the last standard unit vector.
If $x\in\R^n$ then $\|x\|^2:=x^\T x$ is the Euclidean norm.
Finally, $l_{2e}^n$ is the space of all sequences $x:\N_0\to\R^n$, which
are tacitly assumed to be extended as $x_{t}=0$ for $t<0$.

We follow the custom in robust control to express a linear system $x_{t+1}=Ax_t+Bu_t$, $y_t=Cx_t+Du_t$ for $t\in\N_0$ as
$$
\mat{c}{x_{t+1}\\y_t}=\mat{cc}{A&B\\C&D}\mat{c}{x_{t}\\u_t}
\te{or}y=\mas{c|c}{A&B\hl C&D}u.
$$
The latter notation is also used to represent
the input-output map defined by the system. Moreover,
the shorthand notation
$$
\mat{c|c}{
\Al      &\Bl \hl
\Cl      &\Dl  }
\stackrel{T}\longrightarrow
\mat{c|c}{
\Ar&\Br\hl
\Cr&\Dr}
$$
expresses that
$T$ is invertible with
$\Ar=T\Al T^{-1}$, $\Br=T\Bl$, $\Cr=\Cl T^{-1}$, $\Dr=\Dl$. Finally,
we abbreviate Kalman's controllability matrix of the pair $(A,B)\in\R^{n\times(n+m)}$
by $\c{K}(A,B):=(B,AB,\ldots,A^{n-1}B).$

\section{Optimization Algorithms as Feedback Systems}\label{SecII}

\subsection{The Underlying Function Class}

In this paper, we work with the class $\c{S}_{m,L}$ of functions $f:\R^d\to\R$ that are $L$-smooth
and $m$-strongly convex for $m>0$ or just convex if $m=0$. Among the various equivalent ways to express these conditions, the following most intuitive ones do not require any a priori assumptions on differentiability.

\definition{\label{Dsec}Let $L>m\geq 0$ and
$q(x)=\frac{1}{2}\|x\|^2$ for $x\in\R^n$.
Then $\c{S}_{m,L}$ is the set of all $f:\R^d\to\R$ such that
$$
f_m:=f-mq\te{and}f^L:=Lq-f\te{are convex.}
$$
Moreover, let $\c{S}_{m,L}^0:=\{f\in\c{S}_{m,L}\mid f(0)=0,\ \nabla f(0)=0\}$.
}

The latter makes sense since it can be shown that any $f\in\c{S}_{m,L}$ is differentiable \cite{Sch23}. For the gradients of $f$, $f_m$, $f^L$ and any $z\in\R^d$, we record the relation
\eql{sec}{
\mat{c}{\nabla f^L(z)\\\nabla f_m(z)}=
\mat{cc}{L\Id&-\Id\\-m\Id&\Id}\mat{c}{z\\\nabla f(z)}.
}

Among the many known inequalities for $f\in \c{S}_{m,L}$,
the one in the following lemma stands out in allowing for a direct construction of integral quadratic constraints.
It is also underlying the proof of \cite[Lemma 8]{LesRec16} and,
if evaluated at finitely many points, identical to the central inequality in \cite[Theorem 4]{TayHen16}.
\ls{}{The proof is reproduced from \cite{Sch23} in Sec.~\ref{AL2}.}

\lemma{\label{Ldi}
Let $f\in\c{S}_{m,L}$. Then the function $V(x):=(L-m)f_m(x)-q(\nabla f_m(x))$ satisfies
\eql{di}{V(u)-V(y)\leq \nabla f_m(u)^\T[\nabla f^L(u)-\nabla f^L(y)]}
for all $u,y\in\R^d$.
If $f\in\c{S}_{m,L}^0$ then $V$ has a global minimum at $0$ with value $0$, i.e.,  $0=V(0)\leq V(x)$ for all $x\in\R^d$.
}

To support the reader's intuition, we note that \r{di} for $m=0$ and $L\to\infty$
boils down to the subgradient inequality for convex functions.
Since $f_0(x)=f(x)$, $\nabla f_0(x)=\nabla f(x)$ and
$\frac{1}{L}\nabla f^L(x)\to x$ for $L\to\infty$, we infer $\frac{1}{L}V(x)\to f(x)$.
After dividing \r{di} by $L$, we indeed obtain for $m=0$ and $L\to\infty$ the inequality
$f(u)-f(y)\leq \nabla f(u)^T(u-y)$.


\subsection{Optimization Algorithms and Systems}

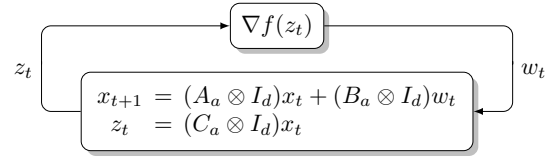
\begin{figure}\center
\scalebox{0.9}{
\begin{tikzpicture}[xscale=1,yscale=1]
\def\dl{11*\dn}
\node[sy3] (g) at (0,0)  {$\arr{ccl}{x_{t+1}&=&(\Aalg\ot \Id)x_t+(\Balg\ot \Id)w_t\\z_t&=&(\Calg\ot \Id)x_t}$};
\node[sy0] (d) at (0,8*\dn)  {$\nabla f(z_t)$};
\draw[<-] (g)-- ($(g) + (2*\dl, 0)$)  |- node[swap]{$w_t$} ($(g) + (2*\dl, .6)$) |- (d);
\draw[->] (g)-- ($(g) + (-2*\dl, 0)$)  |- node[]{$z_t$} ($(g) + (-2*\dl, .6)$) |- (d);
\end{tikzpicture}}
\caption{Feedback representation of optimization algorithms.}\label{fig1}
\end{figure}

For $f\in\c{S}_{m,L}$, we recall that the optimization problem
\eql{opt}{
\inf_{z\in\R^d} f(z)
}
does admit a unique solution $z_*$  \cite{Bec17}. It is also well-known that the gradient descent algorithm
\eql{gd}{
z_{t+1}=z_t-\al\nabla f(z_t)
}
for $\al=\frac{2}{m+L}$ generates a sequence with $\lim_{t\to\infty} z_t=z_*$. We denote the iteration index by ``$t$'' since we want to view \r{gd} as a discrete-time dynamical system for $t$ on the time axis $\N_0$. Even more, \r{gd} can be viewed as the feedback interconnection of the LTI system
$x_{t+1}=x_t-\al w_t$, $z_t=x_t$
with the static nonlinearity
\eql{nl}{
w_t=\nabla f(z_t)
}
for $t\in\N_0$, where $x_t\in\R^d$, $w_t\in\R^d$ and $z_t\in\R^d$
are the state, the input and the output of the linear system. This linear system
can actually be expressed as
\eql{sy}{
x_{t+1}=(\Aalg\otimes I_d)x_t+(\Balg\ot I_d)w_t,\ z_t=(\Calg\ot I_d)x_t
}
where $\Aalg=1$, $\Balg=-\al$ and $\Calg=1$. Here $\ot$ denotes the Kronecker product, which is convenient
to compactly describe general algorithms in the sequel.
In control, the feedback interconnection \r{nl}-\r{sy} is a so-called Lur'e system.
A block diagram of this interconnection is depicted in Fig.~\ref{fig1}.

Accelerated versions of gradient descent include a so-called momentum term. A prominent example
is Nesterov's algorithm with a description
$$
v_{t+2}=v_{t+1}+\be(v_{t+1}-v_t)-\al \nabla f(v_{t+1}+\ga(v_{t+1}-v_{t}))
$$
for suitable real parameters $\al,\be$ and $\ga=\be$ \cite{Nes18}, or the triple momentum version with $\ga\neq\be$ \cite{ScoFre18}. This is nothing but
$$
v_{t+2}=v_{t+1}+\be(v_{t+1}-v_t)-\al w_t,\ z_t=v_{t+1}+\ga(v_{t+1}-v_{t})
$$
in feedback with \r{nl}. Moreover, the latter second order system can be routinely translated into the
first-order description \r{sy} with state $x_t=\col(v_{t+1},v_{t})$ and the matrices
\eql{nes}{
\mat{cc}{\Aalg&\Balg\\\Calg&0}=\mat{cc|c}{1+\be&-\be &-\al \\1&0&0\hl 1+\ga&-\ga &0}.
}
Hence, also Nesterov's recursion can be expressed as \r{nl}-\r{sy}.

These observations provide a strong motivation for investigating the stability properties of the feedback interconnection \r{nl}-\r{sy} for general matrices $(\Aalg,\Balg,\Calg)$.

\subsection{Minimal Convergence Requirement}

In view of the goal to solve \r{opt}, it is a minimal requirement that,
for any $f\in\c{S}_{m,L}$ and any initial condition, the signal $z_t$ of the interconnection \r{nl}-\r{sy} should converge to a limit which satisfies the following first order necessary and sufficient condition for optimality:
\eql{conver}{
\lim_{t\to\infty}z_t=z_*\te{with}\nabla f(z_*)=0.
}

By \r{nl}, this implies $w_*:=\lim_{t\to\infty}w_t=0$.
If $(\Aalg,\Calg)$ is detectable, we infer $\lim_{t\to \infty}x_t=x_*$ and the limit
$(x_*,w_*,z_*)$ satisfies the equilibrium equations
\eql{eq}{x_*=(\Aalg\ot \Id)x_*,\ \ z_*=(\Calg\ot \Id)x_*,\ \ w_*=0.}
If we pick $f(z)=\frac{1}{2}m\|z-z_*\|^2$ for $z_*\in\R^d$ with $z_*\neq 0$, the corresponding
solution $x_*$ of \r{eq} does not vanish, which in turn  shows that $1$ is an eigenvalue of $\Aalg$.
For example in Nesterov's algorithm, this is indeed true since the elements in each row of $\Aalg$ in \r{nes} sum up to one.
As a result, the minimal requirement enforces structural constraints on the algorithm parameters
$(\Aalg,\Balg,\Calg)$.

In general, we argue in \cite[Section 2.2]{SchEbe21} that $(\Aalg,\Calg)$ can be assumed to be detectable without loss of generality. Then the minimal requirement implies that
the algorithm parameters must admit, after a possible state-coordinate change, the structure
\eql{im0}{
\mat{cc}{\Aalg&\Balg\\\Calg&0}=
\mat{cc|c}{\Aalp&\Balp&0\\0&1&1\hl \Calp&\Dalp&0}\text{\ with\ }\det\mat{cc}{\Aalp\!-\!I&\Balp\\ \Calp&\Dalp}\neq 0.
}
The first relation means that the system described with $(\Aalg,\Balg,\Calg)$ is the  series interconnection of

\eql{con}{
\mat{c}{x^c_{t+1}\\z_t}=\mat{cc}{\Aco\ots\Id&\Bco\ots\Id\\\Cco\ots\Id&\Dco\ots\Id}\mat{c}{x^c_{t}\\y_t}
}
and the discrete-time integrator
\eql{pla}{
\mat{c}{x^s_{t+1}\\y_t}=\mat{cc}{\Id&\Id\\\Id&0}\mat{c}{x^s_{t}\\w_t}
}
with the transfer matrix $\frac{1}{\z-1}\Id$. The second condition in \r{im0}
expresses the fact that the pole $\z=1$ of the integrator in the corresponding product of transfer matrices
is not canceled.

In other words, the algorithm's parameters must contain a model of the integrator.
Although not surprising from the perspective of control,
this fact has only been recently clearly emphasized in \cite{MicSch21,SchEbe21} in the realm of algorithm analysis.

As an illustration, for \r{nes} we note that
$$
\mat{cc|c}{1\!+\!\be&-\be &-\al \\1&0&0\hl 1\!+\!\ga&-\ga &0}
\stackrel{\arc\smat{cc}{\ \ 0&1\\-\frac{1}{\al}&\frac{\be}{\al}}}{\longrightarrow}
\mat{cc|c}{\be&-\al &0 \\0&1&1\hl\be(1\!+\!\ga)\!-\!\ga&-\al(1\!+\!\ga)&0}.
$$

Let us now pinpoint the two essential consequences
in case that $(\Aalg,\Balg,\Calg)$ does indeed have the structure \r{im0}:
\enu{
\item If \r{conver} is satisfied for \r{nl}-\r{sy} and all $f\in\c{S}_{m,L}^0$, then
\r{conver} holds for \r{nl}-\r{sy} and all $f\in\c{S}_{m,L}$.
\item The interconnection \r{nl}-\r{sy} is a controlled uncertain system as familiar in robust control.
}

To see 1), we assign to any $f\in\c{S}_{m,L}$ the function $f_*$, defined with the unique $z_*\in\R^d$ satisfying
$\nabla f(z_*)=0$ as
$$
f_*(z):=f(z+z_*)-f(z_*)\te{for}z\in\R^d.
$$
Since $f_*(0)=0$ and $\nabla f_*(0)=0$, we note that $f_*\in\c{S}_{m,L}^0$.
If $x_*$ denotes the unique solution of \r{eq}, the system \r{sy}
can be equivalently transformed into
$$\arr{ccl}{x_{t+1}-x_*&=&(\Aalg\otimes I_d)(x_t-x_*)+(\Balg\ot I_d)w_t,\\ z_t-z_*&=&(\Calg\ot I_d)(x_t-x_*).}$$
To be precise,  the trajectories $(x,w,z)$ of the interconnection \r{nl}-\r{sy} are in one-to-one correspondence
via $\tilde x=x-x_*$, $\tilde w=w$, $\tilde z=z-z_*$  with the trajectories
$(\tilde x,\tilde w,\tilde z)$ of \r{sy} in feedback with
\eql{nlt}{w_t=\nabla f_*(z_t).}
This proves 1).
Even stronger, it shows that exponential stability of
the equilibrium $(\t x_*,\t w_*,\t z_*)=(0,0,0)$
of the interconnection \r{sy}, \r{nlt} is equivalent to
exponential stability of $(x_*,0,z_*)$ of the loop \r{nl}-\r{sy}.

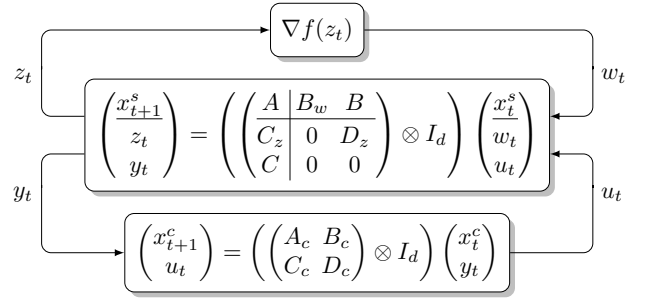
\begin{figure}\center
\scalebox{0.9}{
\begin{tikzpicture}[xscale=1,yscale=1]
\def\dl{2*\dn}
\def\ds{4*\dn}
\node[sy3] (g) at (0,0)  {$\mat{c}{x^s_{t+1}\hl z_t\\y_t}=
\left(
\mat{c|ccc}{\As&\Bsw&B\hl \Csz&0&\Dsz\\\Cs&0&0}
\ot I_d\right)
\mat{c}{x_t^s\hl w_t\\u_t}$};
\tio{i1}{g}{east}{2/3};
\tio{o1}{g}{west}{2/3};
\tio{i2}{g}{east}{1/3};
\tio{o2}{g}{west}{1/3};
\node[sy3,below=\dl of g] (k) {$\mat{c}{x_{t+1}^c\\u_t}=
\left(\mat{cc}{\Aco&\Bco\\\Cco&\Dco}\ot I_d\right)
\mat{c}{x_t^c\\y_t}$};
\node[sy3,above=\dl of g] (d)  {$\nabla f(z_t)$};
\draw[->] (o1)-- ($(o1) + (-\ds, 0)$)  |- node[]{$z_t$} ($(o1) + (-\ds, .6)$) |- (d);
\draw[<-] (i1)-- ($(i1) + (\ds, 0)$)   |- node[swap]{$w_t$} ($(i1) + (\ds, .6)$) |- (d);
\draw[<-] (i2)-- ($(i2) + (\ds, 0)$)   |- node[swap]{$u_t$} ($(i2) + (\ds, -.6)$) |- (k);
\draw[->] (o2)-- ($(o2) + (-\ds, 0)$)  |- node[]{$y_t$} ($(o2) + (-\ds, -.6)$) |- (k);
\end{tikzpicture}
}
\caption{General algorithm synthesis configuration.}\label{fig2}
\end{figure}

Property 2) is seen by redrawing Fig.~\ref{fig1} as in Fig.~\ref{fig2} with
\eql{algpla}{
\mat{c|ccc}{\As&\Bsw&B\hl \Csz&0&\Dsz\\\Cs&0&0}=
\mat{c|cc}{1&1&0\hl 0&0&1\\1&0&0},}
by recalling \r{con}-\r{pla} and using the auxiliary signal $u_t:=z_t$.
Then we indeed recognize the uncertainty $\nabla f$ in the class $\nabla \c{S}_{m,L}^0$, the to-be-controlled plant defined with \r{algpla} and the controller \r{con}.
We also emphasize the simplicity of this plant in the realm of algorithms!
Still, it is relevant to stress  that all our subsequent analysis and synthesis results
even apply to general LTI plants as in Fig.~\ref{fig2}.

\newcommand{\Aalgrho}{\bar A_a}
\newcommand{\Balgrho}{\bar B_a}

\newcommand{\Arho}{\bar A}
\newcommand{\Brho}{\bar B}

\section{Exponential Stability and Passivity}\label{SecIII}

Let us now turn to the development of a test which ensures that the loop \r{nl}-\r{sy} is exponentially stable
for all $f\in\c{S}_{m,L}^0$.

One way is based on the exponential signal weighting map
\eql{trafo}{
T_{\rho^{-1}}(z_0,z_1,z_2,\ldots)=(z_0,\rho^{-1}z_1,\rho^{-2}z_2,\ldots)
}
for some $\rho\in(0,1]$ \cite{DesVid75}. Clearly, $T_{\rho^{-1}}$ is linear and invertible with $T_{\rho^{-1}}^{-1}=T_{\rho}$.
It is then easily checked that the set of trajectories $(x,w,z)$ of \r{sy} are in one-to-one correspondence with trajectories $(\bar x,\bar w,\bar z)$ of the system
\eql{sysrho}{
\arr{ccl}{\bar x_{t+1}&=&(\Aalgrho\ot I_d)\bar x_t+(\Balgrho\ot I_d)\bar w_t,\\
\bar z_t&=&(\Calg\ot I_d)\bar x_t}
}
under the signal transformations $\bar x=T_{\rho^{-1}} x$, $\bar w=T_{\rho^{-1}}w$, and $\bar z=T_{\rho^{-1}}z$
and with the $\rho$-scaled matrices
$
\mat{cc}{\Aalgrho &\Balgrho}:=\rho^{-1}\mat{cc}{A&B}.
$
Similarly, \r{nl} translates into
\eql{nlrho}{
\bar w_t=\bar F(t,\bar z_t)
}
with the static time-varying map $\bar F$ associated to $\nabla f$ through 
\eql{Delrho}{\bar F(t,z):=\rho^{-t}\nabla f(\rho^tz)\te{for}t\in\N_0,\  z\in\R^d.}
The bars should remind us of the fact that that $(\Aalgrho,\Balgrho)$ and $\bar F$ depend on $\rho$.

As a consequence, if the transformed loop is  Lyapunov stable in the sense of $\|\bar x_t\|\leq K\|\bar x_0\|$ for all $t\in\N_0$,
one can conclude that the original loop is  exponentially stable with rate $\rho$ in the sense of
$\|x_t\|\leq K\rho^{-t}\|x_0\|$ for all $t\in\N_0$. This motivates to develop a robust  stability test for the transformed interconnection \r{sysrho}-\r{nlrho}.

We start by deriving what is called a valid integral quadratic constraint (IQC) for the nonlinearity
by  exploiting Lemma~\ref{Ldi}. If $f\in\c{S}_{m,L}^0$, we conclude from \r{di} for $y=0$ that
\eql{dis}{
V(u)\leq \nabla f_m(u)^\T\nabla f^L(u)\te{for all}u\in\R^d.
}
Since $V$ is nonnegative and $\rho>0$, this trivially implies
$$
0\leq \rho^{-t}\nabla f_m(\rho^tz_t)^\T(\rho^{-t}\nabla f^L(\rho^tz_t))
$$
for all $z\in l_{2e}^d$ and $t\in\N_0$. Summation leads to the IQC
\eql{iqc0}{
0\leq \sum_{t=0}^{T-1}\bar F_m(t,\bar z_t)^\T \bar F^L(t,\bar z_t)\te{for all}T\in\N
}
and all sequences $\bar z\in l_{2e}^d$, where $\bar F_m$ and $\bar F^L$ are also defined according  to \r{Delrho}.
Note that the misnomer ``IQC'' results from a similar concept for continuous-time systems, in which
summation is replaced by integration \cite{MegRan97}.
With \r{sec} we infer
\eql{gra}{
\mat{c}{\bar F^L(t,\bar z_t)\\\bar F_m(t,\bar z_t)}=
\mat{cc}{L\Id&-\Id\\-m\Id&\Id}
\mat{c}{\bar z_t\\\bar F(t,\bar z_t)}.
}
This motivates to introduce the static filter
\eql{psi}{
\mat{c}{\bar p_t\\\bar q_t}=
\mat{cc}{L\Id&-\Id\\-m\Id&\Id}
\mat{c}{\bar z_t\\\bar w_t}.
}
If filtering the input-output signals of the nonlinearity \r{nlrho} accordingly,
\r{gra} shows $\bar p_t=\bar F^L(t,\bar z_t)$ and $\bar q_t=\bar F_m(t,\bar z_t)$.
Then \r{iqc0} reads as $\sum_{t=0}^{T-1}\bar q_t^\T \bar p_t\geq 0$ for all $T\in \N$
and can be interpreted as a passivity property \cite{DesVid75} for the outputs of \r{psi}
driven by the signals in \r{nlrho}.

\begin{figure}\center
\scalebox{0.9}{
\begin{tikzpicture}[xscale=1,yscale=1]
\def\dl{11*\dn}
\node[sy3] (g) at (0,0)  {$\arr{ccl}{\bar x_{t+1}&=&(\Aalgrho\ot \Id)\bar x_t+(\Balgrho\ot \Id)\bar w_t\\\bar z_t&=&(\Calg\ot \Id)\bar x_t}$};
\node[sy0] (d) at (-12*\dn,8*\dn)  {$\bar F(t,\bar z_t)$};
\node[sys={5*\dn}{6*\dn}] (f) at (8*\dn,15*\dn) {
$\arr{ccc}{\bar p_t&=&L\bar z_t-\bar w_t\\\bar q_t&=&-m\bar z_t+\bar w_t}$};
\tio{i1}{f}{west}{1/3};
\tio{i2}{f}{west}{2/3};
\tio{o1}{f}{east}{1/3};
\tio{o2}{f}{east}{2/3};


\draw[<-] (i1) -| ([xshift=3*\dn] d.east);
\draw[<-] (i2) -| ([xshift=-3*\dn] d.west);

\draw[->] (o2)--node[]{$\bar p_t$} ([xshift=5*\dn] o2) ;
\draw[->] (o1)--node[swap]{$\bar q_t$} ([xshift=5*\dn] o1) ;

\draw[<-] (g)-- ($(g) + (2*\dl, 0)$)  |- node[swap]{$\bar w_t$} ($(g) + (2*\dl, .6)$) |- (d);
\draw[->] (g)-- ($(g) + (-2*\dl, 0)$)  |- node[]{$\bar z_t$} ($(g) + (-2*\dl, .6)$) |- (d);
\end{tikzpicture}}
\caption{Block diagram of transformed loop \r{sysrho}-\r{nlrho} with signals filtered according to \r{psi}.}\label{fig3a}
\end{figure}
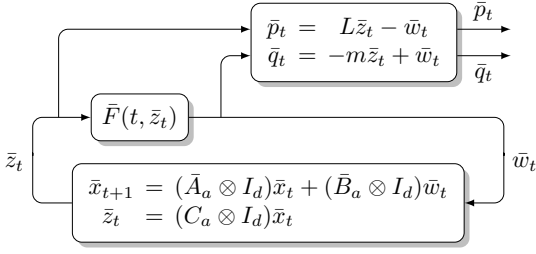

In view of Fig.~\ref{fig3a}
and motivated by the passivity theorem,
we  expect that stability of \r{sysrho}-\r{nlrho} is guaranteed
in case that $\sum_{t=0}^{T-1}\bar q_t^\T \bar p_t<0$ holds for all $T\in \N$
along the input-output trajectories of the linear system \r{sysrho}
filtered with \r{psi}.

To make this precise, we start by emphasizing that the stability test itself is
formulated for $d=1$, while the conclusions are drawn for arbitrary dimensions
$d\in\N$. Throughout the paper we slightly abuse the notation and do not indicate
the dependence of system signals on $d$.

For $d=1$, we note that the input-output signals of \r{sysrho}
filtered with \r{psi} satisfy
$\bar w_t=m\Calg\bar x_t+\bar q_t$ and $\bar p_t=L\Calg\bar x_t-\bar w_t$.
With an identical state-trajectory, we hence infer
\eql{sya}{
\mat{c}{\bar x_{t+1}\\\bar p_t}=
\underbrace{\mat{cc}{\Aalgrho+\Balgrho m\Calg&\Balgrho\\\ (L-m)\Calg&-1}}_{\footnotesize\mat{cc}{\Acl&\Bcl\\ \Ccla&\Dcla}}
\mat{c}{\bar x_{t}\\ \bar q_t}.
}

As a consequence, also for a general $d\in\N$,
the trajectories of \r{sysrho} filtered with \r{psi}
satisfy
\eql{sysrhof}{
\mat{c}{\bar x_{t+1}\\ \bar p_t}=
\mat{cc}{\Acl\otimes \Id&\Bcl\otimes \Id\\ \Ccla\otimes \Id&\Dcla\otimes \Id}\mat{c}{\bar x_t\\ \bar q_t}.
}
This leads to our first analysis result, which involves a passivity property of the system \r{sya} and hence also of \r{sysrhof}.

\theorem{\label{Trs}
Suppose there exists some $\Xcl=\Xcl^\T$ with
\mul{\label{lmia}
\Xcl\cg0\text{\ and\ }
\mat{cc}{\Acl&\Bcl\\I&0}^\T\!
\mat{cc}{\Xcl&0\\0&-\Xcl}
\mat{cc}{\Acl&\Bcl\\I&0}+\\+
\mat{cc}{\Ccla&\Dcla\\0&1}^\T\!
\mat{cc}{0&1\\1&0}
\mat{cc}{\Ccla&\Dcla\\0&1}
\cl 0.
}
Then there exists a constant $K$ such that, for any $f\in\c{S}_{m,L}^0$,
all trajectories of the original loop \r{nl}-\r{sy} satisfy
\eql{expsta}{
\|x_t\|\leq K\rho^{t}\|x_0\|\te{for all}t\in\N_0.
}
For $\rho=1$, it is also assured that $\lim_{t\to\infty}x_t=0$ holds true.}
\ls{}{The dissipativity-based proof is found in Sec.~\ref{ATrs}.}

Before addressing the practical application of this robust stability test, we discuss
how to substantially improve it by the incorporation of so-called stability multipliers.

\section{Dynamic Integral Quadratic Constraints}\label{SecIV}

It is a classical idea \cite{WilBro68,Zamfal68}
to improve Theorem~\ref{Trs} by imposing a passivity condition after filtering the signal $\bar p$ in \r{psi} with a
causal and stable time-invariant system. In this context, such a filter is often called a stability multiplier \cite{DesVid75}.

In fact, passing the signal $\bar p_t=\bar F^L(t,\bar z_{t})$ through a delay of time $\nu\in\N$ leads to
$\bar r_t=\bar F^L(t-\nu,\bar z_{t-\nu})$ (where we recall our convention that $\bar z_{t-\nu}=0$ and hence $\bar r_{t-\nu}=0$ for $t<\nu$.) The following IQC incorporates this delayed signal $\bar r_t$ and is, again,
a rather immediate consequence of Lemma~\ref{Ldi}.
\lemma{\label{Lzf}Let $f\in\c{S}_{m,L}^0$, $\rho\in(0,1]$ and $\nu\in\N$. Then
\eql{iqc1}{
0\leq\sum_{t=0}^{T-1}\bar F_m(t,\bar z_t)^\T(\bar F^L(t,\bar z_t)-\rho^{\nu}\bar F^L(t-\nu,\bar z_{t-\nu}))
}
holds for all $T\in\N$ and all signals $\bar z\in l_{2e}^d$.
}

\ls{}{The proof is found in Sec.~\ref{AL4}.}
A conic combination of \r{iqc0} and \r{iqc1} for $\nu\in\N$ leads to the
IQC with more general filters in the following lemma
\ls{\!\!.}{ as proved in Sec.~\ref{AL5}.}

\lemma{\label{Lzfd}
Let $f\in\c{S}_{m,L}^0$, $\rho\in(0,1]$ and suppose that $\la_0,\la_1,\ldots\in\R$ satisfy
\eql{iirc}{
\la_\nu\leq 0\te{for all}\nu\in\N\te{and}\sum_{\nu=0}^\infty\rho^{-\nu}\la_\nu>0.
}
For $\bar z\in l_{2e}^d$, let $\bar p_t=\bar F^L(t,\bar z_t)$ be passed through the filter
\eql{iir}{
\bar r_t=\sum_{\nu=0}^{t} \la_{\nu} \bar p_{t-\nu}\te{for}t\in\N_0.
}
With $\bar q_t:=\bar F_m(t,\bar z_t)$,
the signals $\bar r,\bar q$ then satisfy the IQC
\eql{did}{
0\leq\sum_{t=0}^{T-1}\bar q_t^\T \bar r_t\te{for all}T\in\N.
}
}

Due to the incorporation of the dynamic filter \r{iir},
we call \r{did} a dynamic IQC. Note that \r{iirc} implies $\la_0>0$.

The infinite impulse response filters \r{iir}
are subject to the infinite number of constraints
\r{iirc}. To overcome this trouble for the purpose of computations, we proceed with
filters that have
a state-space realization $(\Af,\Bf,\Cfla,\Dfla)$
with a fixed pole-pair $(\Af,\Bf)\in\R^{l\times(l+1)}$
and free filter coefficients collected in $(\Cfla,\Dfla)\in\R^{1\times l}\times\R$
such that
\eql{fil}{\la_0=\Dfla\te{and}\la_{\nu+1}=\Cfla\Af^\nu\Bf\text{\ for\ }\nu\in\N_0.}
The boldface notation reminds us of the fact that
$(\Cfla,\Dfla)$ is a decision variable in
the subsequent stability test.
For a suitable choice of $(\Af,\Bf)$, 
we now establish that
the infinitely many constraints \r{iirc} on the Markov parameters
\r{fil} can be expressed by a finite number of linear ones.

\lemma{\label{Lpf}
Fix $(\Af,\Bf)=(C_\al,e_l)$ with a real polynomial
$\al(\z)=\alpha_0+\alpha_1\z+\cdots+\alpha_{l-1}\z^{l-1}+\z^l$
of degree $l\in\N$ having all its roots in
$\D_\rho:=\{\z\in\C\mid |\z|<\rho\}$ and with coefficients satisfying
$\al_0,\ldots,\al_{l-1}\leq 0.$

Then, for any pair $(\Cfla,\Dfla)\in\R^{1\times l}\times\R$, the Markov parameters \r{fil} satisfy the constraints \r{iirc} iff
\eql{lac}{
\Cfla\c{K}(\Af,\Bf)\leq 0\te{and}
\Dfla+\Cfla(\rho I-\Af)^{-1}\Bf>0.}
Moreover, \r{lac} implies that $\Dfla>0$ and that all eigenvalues of $\Af-\Bf\Dfla^{-1}\Cfla$ are located in $\D_\rho$.
}

\ls{}{The proof is given in Sec.~\ref{AL6}.}

\begin{figure}\center
\scalebox{0.9}{
\begin{tikzpicture}[xscale=1,yscale=1]
\def\dl{4*\dn}
\node[sy0] (d) at (0,0)  {$\mat{c}{\bar z \\\bar F(.,\bar z)}$};
\tio{od1}{d}{east}{1/3};
\tio{od2}{d}{east}{2/3};
\node[sys={5*\dn}{6*\dn},right=\dl of d] (f)  {
$\mat{cc}{L\Id&-\Id\\-m\Id&\Id}$};
\tio{i1}{f}{west}{1/3};
\tio{i2}{f}{west}{2/3};
\tio{o1}{f}{east}{1/3};
\tio{o2}{f}{east}{2/3};

\node[sy0,right= 1.5*\dl of o2] (iir) {
\r{fild}
};

\draw[<-] (i1) --node[]{$\bar w$} (od1);
\draw[<-] (i2) --node[swap]{$\bar z$} (od2);

\draw[->] (o2)--node[pos=.35]{$\bar p$} (iir);
\draw[->] (o1)--node[swap]{$\bar q$} ([xshift=\dl] o1) ;

\draw[<-] (d.west)--node[swap]{$\bar z$} ([xshift=-\dl] d.west) ;

\draw[->] (iir.east)--node[]{$\bar r$} ([xshift=\dl] iir.east) ;
\end{tikzpicture}}
\caption{Signals in dynamic IQC of nonlinearity.}\label{fig3b}
\end{figure}
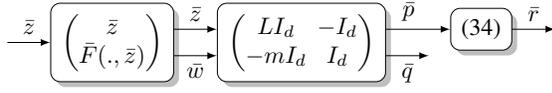

To summarize, for a fixed polynomial $\al$ of degree $l$ as in Lemma~\ref{Lpf},
we work from now on with the filter matrices
\eql{filp}{
\Af:=C_\al,\ \ \Bf:=e_l,\ \ \Cfla\in\R^{1\times l},\ \ \Dfla\in\R
}
such that $(\Cfla,\Dfla)$ satisfies the constraints \r{lac}.
It is then assured that the trajectories of \r{nlrho} filtered by \r{psi} and
\eql{fild}{
\mat{c}{\xi_{t+1}\\\bar r_t}=\mat{cc}{\Af\ots\Id&\Bf\ots\Id\\\Cfla\ots\Id&\Dfla\ots\Id}\mat{c}{\xi_{t}\\\bar p_t},\ \ \xi_0=0
}
(see Fig.~\ref{fig3b}) satisfy the passivity condition \r{did} for all $f\in\c{S}_{m,L}^0$.
In this way, we have identified a whole nicely parameterized convex family of valid dynamic IQCs
for the nonlinearity \r{nlrho} involving the filters or multipliers \r{fild}.


\section{Robust Stability Analysis with Dynamic IQCs}
\label{SecV}

\begin{figure}\center
\scalebox{0.9}{
\begin{tikzpicture}[xscale=1,yscale=1]
\def\dl{4*\dn}
\node[sy3] (d) at (0,0)  {$\mas{c|c}{\Aalgrho&\Balgrho\hl \Calg&0\\0&1}$};
\tio{od1}{d}{east}{1.53/4};
\tio{od2}{d}{east}{2.47/4};
\node[sys={5*\dn}{6*\dn},right=\dl of d] (f)  {
$\mat{cc}{L&-1\\-m&1}$};
\tio{i1}{f}{west}{1/3};
\tio{i2}{f}{west}{2/3};
\tio{o1}{f}{east}{1/3};
\tio{o2}{f}{east}{2/3};

\node[sy0,right= 1.5*\dl of o2] (iir) {$\mas{c|c}{\Af&\Bf\hl \Cfla&\Dfla}$};

\draw[<-] (i1) -- node[]{$\bar w$} (od1);
\draw[<-] (i2) -- node[swap]{$\bar z$} (od2);

\draw[->] (o2)--node[pos=.35]{$\bar p$} (iir);
\draw[->] (o1)--node[swap]{$\bar q$} ([xshift=\dl] o1) ;

\draw[<-] (d.west)--node[swap]{$\bar w$} ([xshift=-\dl] d.west) ;

\draw[->] (iir.east)--node[]{$\bar r$} ([xshift=\dl] iir.east) ;
\end{tikzpicture}}

%
%
%
%
%

\caption{Signals in filtered linear system.}\label{fig3c}
\end{figure}
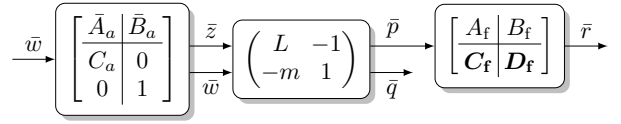

In alignment with Sec.~\ref{SecIII}, the robust stability test
with dynamic IQCs is now formulated for the correspondingly filtered linear system as depicted
in Fig.~\ref{fig3c}. If recalling \r{sya}, we are lead to the system
$$
\bar r=
\mas{c|c}{\Af&\Bf\hl \Cfla&\Dfla}
\mas{c|c}{\Aalgrho+\Balgrho m\Calg&\Balgrho\hl (L-m)\Calg&-1}\bar q.
$$
A realization of this series interconnection is given  as
\eql{sycf}{\arc
\mat{c|c}{\Acl&\Bcl\hl \Ccla&\Dcla}:=
\mat{cc|ccc}{
\Af&\Bf (L-m)\Calg&-\Bf  \\
0&\rho^{-1}(\Aalg +\Balg m\Calg)&\rho^{-1}\Balg\hl
\Cfla& \Dfla(L-m)\Calg&-\Dfla}.
}

\theorem{\label{Trs2}
Fix some $\rho\in(0,1]$ and $(\Af,\Bf)$ with the properties in Lemma~\ref{Lpf}.
If there exist an $\Xcl=\Xcl^\T$ and a
filter parameter $(\Cfla,\Dfla)$ with \r{lac}
such that \r{lmia} is satisfied,
then the same conclusions can be drawn as in Theorem~\ref{Trs}.
}

\ls{}{The proof is given in Section~\ref{ATrs2}.
Moreover, some explanatory connections of Theorems~\ref{Trs} and \ref{Trs2} with
classical passivity-based stability tests
are discussed in Sec.~\ref{A0}.}

We use boldface letters to highlight that
$(\Ccla\ \Dcla)$ depends affinely on the decision variables $(\Cfla,\Dfla)$. Hence
\r{lmia} constitutes a genuine LMI in the variables
$(\Cfla,\Dfla)$ and $\Xcl$.

We emphasize that Theorem~\ref{Trs2} encompasses Theorem~\ref{Trs} for $l=0$, which
means that $\Af,\Bf,\Cfla$ are empty matrices; since the LMI \r{lmia} is homogeneous in $(\Xcl,\Dcla)$,
we can indeed fix $\Dcla=\Dfla$ to the value $1$  without loss of generality.

After picking the targeted convergence rate $\rho\in(0,1]$ and either $l=0$ or
a characteristic filter polynomial $\al$ of degree $l\in\N$ as in Lemma~\ref{Lpf},
feasibility of the LMIs \r{lac} and \r{lmia} (defined with the matrices \r{sycf}) thus guarantees
exponential stability of \r{nl}-\r{sy} with rate $\rho$
for any $f\in\c{S}_{m,L}^0$ and any dimension $d\in\N$.
A more detailed practical recipe for how to apply Theorem~\ref{Trs2} can be extracted from
Sec.~\ref{Ssyn}.

Note that the dimension of the LMI \r{lmia} is independent from $d\in\N$. Moreover, it is
remarkable that even only a few filter states ($l=1,2,3$) can substantially improve the test over $l=0$.
For $\al(\z)=\z^l$ with $l=1$, Theorem~\ref{Trs2} encompasses the algorithm analysis approach in \cite{LesRec16}.

Related IQC results to ensure robust exponential loop stability can be found in \cite{BocLes15,HuSei16}, while
a more general dissipativity-based robustness analysis framework is exposed in the recent survey article \cite{Sch22} and the references therein.

\section{Convex Algorithm Synthesis}\label{Ssyn}

Let us now turn to the design problem for the
interconnection in Fig.~\ref{fig2} with the general plant
\eql{gpla}{
\mat{c}{z\\y}=
\mas{c|ccc}{\As&\Bsw&B\hl \Csz&0&\Dsz\\\Cs&0&0}
\mat{c}{w\\u}
}
where $\As\in\R^{n\times n}$. Again, we slightly abuse notation by not indicating the dependence
of the signals on $d$.
As for analysis, we first construct the modified system description
in order to formulate the key synthesis result.
To this end, we transform all signals of \r{gpla} according to \r{trafo} to get 
$$\arc
\mat{c}{\bar z\\\bar y}=
\mas{c|ccc}{\rho^{-1}\As &\rho^{-1}\Bsw&\rho^{-1}\Bs\hl \Csz&0&\Dsz\\\Cs&0&0}
\mat{c}{\bar w\\\bar u}
.$$
Filtering the signal $\col(\bar z,\bar w)$ as in \r{psi} with $d=1$ leads to
$$
\arc
\mat{c}{\bar p\\\bar y}\!=\!\!
\underbrace{\mas{c|ccc}{
\rho^{-1}\!(\As\!+\!\!\Bsw m\Csz)&\rho^{-1}\!\Bsw&\rho^{-1}\!(\Bs\!+\!\!\Bsw m\Dsz)\hl
(L-m)\Csz&-1 & (L-m)\Dsz\\C&0&0}}_{\footnotesize\mas{c|ccc}{\At&\Bt_w&\Bt\hl \Ct_z&\Dt_{zw}&\Et\\\Ct&0&0}}\!\!
\mat{c}{\bar q\\\bar u}\!.
$$
This plant filtered with
\r{fild} for $d=1$  reads as
\eql{gpf}{
\mat{c}{\bar r\\\bar y}=
\underbrace{
\mas{cc|ccc}{\Af&\Bf \Ct_z&\Bf\Dt_{zw} &\Bf \Et\\0&\At&\Bt_w&\Bt\hl\Cfla& \Dfla\Ct_z&\Dfla\Dt_{zw} &\Dfla \Et\\0&C&0&0}}_{
\footnotesize
\mas{c|cc}{
\Ah   &\Bh_w &\Bh\hl
\Chla &\Dhla &\Ehla\\
\Ch   &0     &0   }
}
\mat{c}{\bar q\\\bar u}.
}
For the system \r{gpf}, we now pick a controller
\eql{coh}{
\bar u=\mas{c|c}{\Ach&\Bch\hl \Cch&\Dch}\bar y.
}
Then the resulting interconnection admits the description
\eql{gclh}{
\bar r=
\mas{cc|c}{
\Ah+\Bh\Dch  \Ch   &\Bh\Cch   &\Bh_w  \\
\Bch \Ch           &\Ach      &0       \hl
\Chla+\Ehla\Dch\Ch &\Ehla\Cch &\Dhla }
\bar q=:
\mas{c|c}{\Acl&\Bcl\hl \Ccla&\Dcla}\bar q,
}
where all the bold matrices depend affinely on
$(\Cfla,\Dfla)$.

If \r{coh} is a controller for which the controlled system \r{gclh} satisfies the hypotheses of Theorem~\ref{Trs2}, it is not difficult to verify that the original plant \r{gpla} controlled with
\eql{corho}{
u=\mas{c|c}{\rho\Ach&\rho\Bch\hl \Cch&\Dch}y
}
(see Figure~\ref{fig2}) satisfies \r{expsta}
for any $x_0$ and any $f\in\c{S}_{m,L}^0$.

Let us now recap a slight variant of a seminal result obtained in
\cite{GahApk94,IwaSke94}, a convex solution for the design problem if $(\Cfla,\Dfla)$
satisfying \r{lac} is {\em held fixed}. To this end, we pick so-called annihilator matrices $\Uh$ and $\Vh$ with
\eql{ann2}{
\Uh=\diag(\Ch_\bot\ 1)
\te{and}
\Vh^\T=\mat{cc}{\Bh^\T&{\bm {\Dh_z}^\T}}_\bot
}
where $M_\bot$ means that the columns of this matrix form a basis of the kernel of the matrix $M$.
Then there exist a controller \r{coh} for \r{gpf} such that the closed loop system \r{gclh}
renders the analysis LMIs in Theorem~\ref{Trs2} feasible
iff there exist symmetric matrices $X$ and $\Yh$ which satisfy
\begin{gather}
\label{lmip}
\left[\bul\right]^\T\!\!\mat{cc|cc}{
X &0 & 0&0\\
0 &-X& 0&0\hl
0&0&0&1\\0&0&1&0}
\left[\mat{cc}{
\Ah&\Bh_1\\
I&0\hl
\Chla&\Dhla\\
0&1}\Uh\right]\cl 0,\\
\label{lmisd}
\left[\Vh\mat{cc|cc}{-I&\Ah&0&\Bh_1\\0& \Chla&-1&\Dhla}\right]
\mat{cc|cc}{\Yh&0&0&0\\0&-\Yh&0&0\hl 0&0&0&1\\0&0&1&0}
\left[\bul\right]^\T\cg 0,\\
\label{lmisc}
\mat{cc}{\Yh&I\\I&X}\cg 0.
\end{gather}
For reasons of space, we use the bullet notation to indicate that one should substitute (on the left/right)
the respective matrix in square brackets (on the right/left) to render the  inequalities symmetric.
For fixed $(\Cfla,\Dfla)$,  these constraints are affine in $X$ and $\Yh$.
However, this nice structural property is destroyed for
\r{lmisd} if viewing  $(\Cfla,\Dfla)$ as an additional decision variable.

\renewcommand{\Al}{\Af}
\renewcommand{\Bl}{\Bf}
\renewcommand{\Cl}{\Cfla}
\renewcommand{\Dl}{\Dfla}

\renewcommand{\Ar}{\At}
\renewcommand{\Br}{\Bt}
\renewcommand{\Cr}{\Ctz}
\renewcommand{\Dr}{\Et}

To overcome this trouble, we note that
the subsystem $\bar u\to\bar r$ of \r{gpf} is actually given by
\eql{subsys}{
\bar r=
\mas{c|ccc}{
\Af  &\Bf\hl
\Cfla&\Dfla}
\mas{c|ccc}{
\At   &\Bt\hl
\Ctz &\Et }\bar u=
\mas{c|ccc}{
\At   &\Bt\hl
\Ctz &\Et }
\mas{c|ccc}{
\Af  &\Bf\hl
\Cfla&\Dfla}
\bar u,
}
a series interconnection of two commuting SISO systems.
As the key to convexification, we exploit the fact that this commutation property is reflected by a
state-coordinate change for the corresponding natural realizations as
$$
\mas{cc|ccc}{
\Al  &\Bl \Cr   &\Bl \Dr     \\
0    & \Ar     &\Br             \hl
\Cl  &\Dl  \Cr  &\Dl \Dr }
\stackrel{\Tla
}{\longrightarrow}
\mas{cc|ccc}{
\Al         &0              &\Bl   \\
\Br\Cl        &\Ar            &\Br\Dl  \hl
\Dr\Cl        &\Cr            &\Dr\Dl  }
$$
with the specifically structured transformation matrix
\eql{Tra}{
\Tla=\mat{cc}{L^{-1}&-L^{-1}K\\ \Nla L^{-1}&\Mla\!-\!\Nla L^{-1}K}=\mat{c}{\Tf\\ \Ttla}.
}
\ls{}{The precise result is formulated in Lemma~\ref{Lcom} in Section~\ref{Scom}.}

An analogous commutation property is used for convexification in \cite{SchEbe21} based on the Youla-Parametrization and in the state-space approach of \cite{HolSch21c}. The latter is confined
to $\al(\z)=\z^l$ for the characteristic polynomial of $\Af$, which induces some limitation on the plant \r{gpf}.
None are required in the next result due to the novel flexibility of choosing $\al$.

\theorem{\label{Tsyn}Pick $\al$ as in Lemma~\ref{Lpf} such that the eigenvalues of $\Af$
are different from the eigenvalues of $\At$ and from the zeros of $\Ct_z(\z I-\At)^{-1}\Bt+\Et$.
With the solutions $K,L,\Mla$ and $\Nla$ of the linear equations
\enu{
\item $\Af K-K\At+\Bf \Ctz=0,$
\item $L\,\c{K}(\Af,\Bf)=\c{K}(\Af,\Bf \Et-K\Bt)$,
\item $\Mla\al(\At)=\Dfla\al(\At)+(\Cfla\ot I_l)\col(I,\At,\ldots,\At^{l-1}),$
\item $\At\Nla-\Nla\Af+\Bt\Cfla=0$,
}
define $\Ttla=\mat{cc}{\Nla L^{-1}&\Mla\!-\!\Nla L^{-1}K}.$
Moreover, let
\eql{ann}{
\Uh=\diag(I_l,\Ct_\bot,1)
\te{and}
V^\T=\mat{cc}{\Bt^\T&\Et^\T}_\bot.
}
Then the following statements are equivalent.
\begin{enumerate}
\item[(a)] There exists a controller \r{coh} for the plant \r{gpf} such that the controlled interconnection \r{gclh}
satisfies the hypotheses in Theorem~\ref{Trs2}.
\item[(b)]
There exist $(\Cfla,\Dfla)$ with  \r{lac} and symmetric matrices $X$, $\Yt$ satisfying the LMIs \r{lmip} and
\eql{lmid}{\arc
\left[V
\mat{cc|cc}{
-I &\At   &0  &\Ttla  \Bh_1\\
0  &\Ctz &-1 &\Dhla       }\right]
\mat{cc|cc}{\Yt&0&0&0\\0&-\Yt&0&0\hl 0&0&0&1\\0&0&1&0}\left[\bul\right]^\T\cg 0,
}
\eql{lmic}{
\mat{ccc}{\Yt&\Ttla \\\Ttla^\T &X}\cg 0.}
\end{enumerate}
}

By its very definition, $\Ttla$ depends {\em affinely} on $(\Cfla,\Dfla)$ and, thus,
the constraints \r{lac}, \r{lmip} and \r{lmid}-\r{lmic} constitute {\em affine constraints on
all decision variable} $(\Cfla,\Dfla)$, $X$ and $\t Y$.
Hence, their feasibility can be verified by standard SDP-solvers.
We emphasize that the complexity of these synthesis LMIs is determined by the dimensions of $A$ and $\Af$ only.

Note that the set of all eigenvalues and zeros of $\At$ and $\Ct_z(\z I-\At)^{-1}\Bt+\Et$, respectively,
are given by $\rho^{-1}\La$ and $\rho^{-1}\Xi$ with two finite and $\rho$-independent sets $\La,\Xi\subset\C$.
For a particular choice of $(\Af,\Bf)$, let us now summarize a concrete procedure for the synthesis of controllers as follows:
\enu{
\item Fix $0<m<L$. If $0\not\in\La\cup\Xi$ set $z_0=0$. Otherwise
choose $z_0>0$ close to zero such that $|\la|>z_0$ holds for all $\la\in(\La\cup\Xi)\setminus\{0\}$.
\item Pick $l\in\N$ and $\Af:=C_\al$, $\Bf:=e_l$ for $\al(\z):=\z^l-z_0^l$.
\item For any $\rho\in(z_0,1]$,  set up the system of LMIs \r{lac}, \r{lmip}, \r{lmid}-\r{lmic} in the variables $(\Cfla,\Dfla)$, $X$ and $\Yt$.
\item By bisection, determine
the best possible (infimal) rate $\rho_*\in[z_0,1]$ such that the resulting LMIs are feasible.
\item For some $\rho\in(\rho_*,1]$ close to the optimal value $\rho_*$, set up the plant
\r{gpf} with $(\Cfla,\Dfla)$ as obtained from a feasible solution of these LMIs.
\item Design a controller \r{coh} such that the LMI \r{lmia} for the closed-loop system \r{gclh} is feasible in $\c{X}$.
\item Define $(\Aco,\Bco,\Cco,\Dco):=(\rho\Ach,\rho\Bch,\Cch,\Dch)$.
}
Since the assumptions on $\al$ in Theorem~\ref{Tsyn} are satisfied,
it is possible to set up the LMIs in Step 3); indeed, all the eigenvalues of $\Af$ have absolute value $z_0$ and, hence, none of them is contained in $\rho^{-1}\La\cup\rho^{-1}\Xi$.
As for Theorem~\ref{Trs2}, the case $l=0$ is covered with empty matrices $\Af$, $\Bf$, $\Cfla$ and $\Dfla=1$, which boils down to choosing $\Ttla=I$ to set up the LMIs in Step 3). Moreover, Theorem~\ref{Tsyn} guarantees that a controller as in Step 6) does indeed
exist (possibly after a slight perturbation of $(\Cfla,\Dfla)$ as seen in the proof in \cite{SchEbe23a}.)

A numerically stable and constructive procedure to design a controller as in Step 6) and based on the classical synthesis conditions \r{lmip}-\r{lmisc} is found, e.g., in \cite{Gah94}.


With the controller in Step 7), for any $f\in\c{S}_{m,L}^0$ and any dimension $d\in\N$,
it is guaranteed that
all trajectories of the interconnection in Fig.~\ref{fig2} decay exponentially with rate $\rho$.



\section{Numerical Illustrations}\label{SNum}

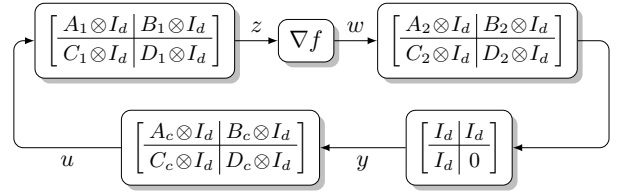
\begin{figure}\center
\scalebox{.9}{
\begin{tikzpicture}[xscale=1,yscale=1]
\def\dl{4*\dn}
\def\fa{3}
\node[sy3] (g) at (0,0)  {\small$\mas{c|c}{\Aalp\ots \Id&\Balp\ots \Id\hl \Calp\ots \Id&\Dalp\ots \Id}$};
\node[sy3,right=2*\dl of g] (k) {\small$\mas{c|c}{\Id&\Id\hl \Id&0}$};
\node[sy0] (d) at (8*\dn,10*\dn)  {$\nabla f$};

\node[sy3,left=1*\dl of d] (g1) {\small$\mas{c|c}{A_1\ots \Id&B_1\ots \Id\hl C_1\ots \Id&D_1\ots \Id}$};
\node[sy3,right=1*\dl of d] (g2) {\small$\mas{c|c}{A_2\ots \Id&B_2\ots \Id\hl C_2\ots \Id&D_2\ots \Id}$};

\node[below =\dn of g] (tg)   {};

\draw[->] (g.west)-- node{$u$} ($(g.west) + (-10*\dn, 0)$)  |-  (g1);
\draw[<-] (k.east)--  ($(k.east) + (9*\dn, 0)$)  |- (g2);
\draw[->] (k) -- node{$y$} node{} (g) ;
\draw[->] (g1) -- node{$z$} (d) ;
\draw[->] (d) -- node{$w$} (g2) ;
\end{tikzpicture}}
\vspace*{-3ex}
\caption{Optimization over communication channels. }\label{fig4}
\end{figure}

We illustrate our results by an extremum control problem.
The purpose is minimize any $f\in\c{S}_{m,L}$ over communication channels.
Concretely, the algorithm needs to transmit
the actual iterate $u_t$ via a channel modeled by an
LTI system with transfer function $G_1(\z)=
C_1(\z I-A_1)^{-1}B_1+D_1$
to generate $z_t$. This is fed into the
gradient to return $w_t=\nabla f(z_t)$. In turn, this signal is communicated
back to the algorithm via a channel
with transfer function $G_2(\z)=
C_2(\z I-A_2)^{-1}B_2+D_2$. To enforce integral action, we are led to
the configuration in Fig.~\ref{fig4} with a to-be-designed
controller $(\Aalp,\Balp,\Calp,\Dalp)$.
To avoid cancelation of the integrator's pole in the loop, we assume that
$\det(A_2-I)\neq 0$ and $G_2(1)\neq 0$.

Since the configuration in Fig.~\ref{fig4} can be subsumed to
the one in Fig.~\ref{fig2}, we can follow the procedure in Sec.~\ref{Ssyn} to
compute optimal rates and close-to-optimal algorithms for this optimization
problem under communication constraints.

If choosing $G_1(\z)=G_2(\z)=1$,
Fig.~\ref{fig4} is identical to Fig.~\ref{fig2} for standard optimization.
With $z_0=0$ in our synthesis procedure, we recover
both the optimal convergence rates and the algorithm parameters
for gradient decent ($l=0$) and $\rho_{\rm tm}:=1-\sqrt{\frac{m}{L}}$ for
the triple momentum algorithm ($l=1$) \cite{Nes83,Nes18,ScoFre18}.
Remarkably, Theorem~\ref{Tsyn} permits to show that
$\rho_{\rm tm}$ is indeed the best possible rate
that is achievable among all algorithms and any $l\in\N_0$ \cite[Corollary 4.8]{SchEbe21}.

\begin{figure}\center
\includegraphics[width=.45\textwidth,trim=0 0 0 0,clip=true]{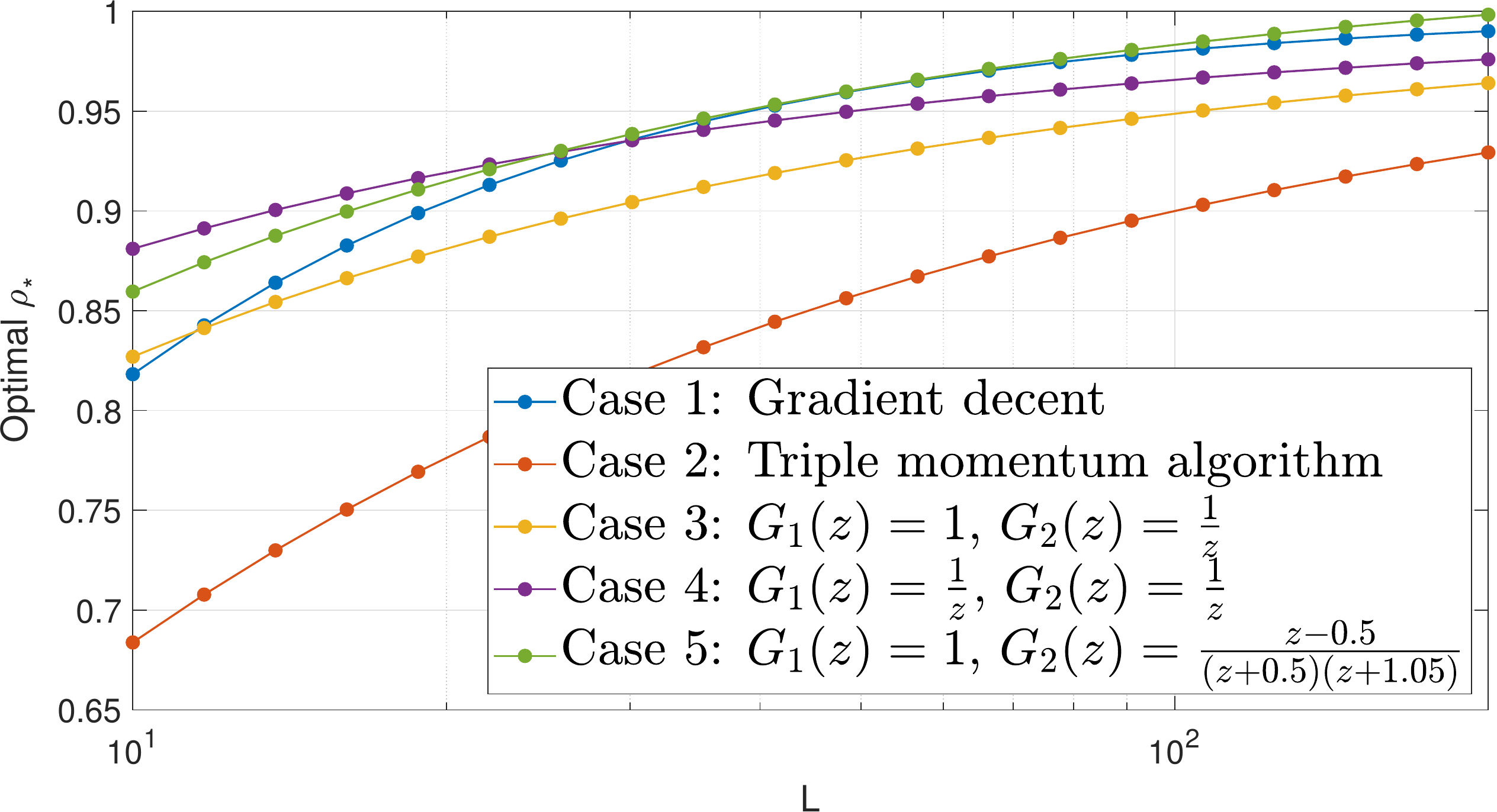}
\caption{Optimal convergence rates plotted over $L$ for $m=1$, $z_0=10^{-2}$ and $l=2$
and for various communication dynamics.}\label{figSim1}
\end{figure}

In the subsequent numerical experiments, we pick $m=1$,
$z_0=10^{-2}$, $l=2$ and compute the optimal rates $\rho_*$
for the communication filters in Fig.~\ref{figSim1}.
The results are plotted over the so-called condition number $L$ of the class $\c{S}_{1,L}$.

If compared to the tripel momentum algorithm (Case 2), the rates increase if the
gradients are processed with a one-step delay (Case 3), but they are still mostly better
than for gradient descent (Case 1). For lager values of $L$, this is even true for delays in both channels (Case 4).
Our approach allows for unstable dynamics in the optimization loop
(Case 5), which affects the convergence rates adversely.

\begin{figure}\center
\includegraphics[width=.48\textwidth,trim=120 15 100 0,clip=true]{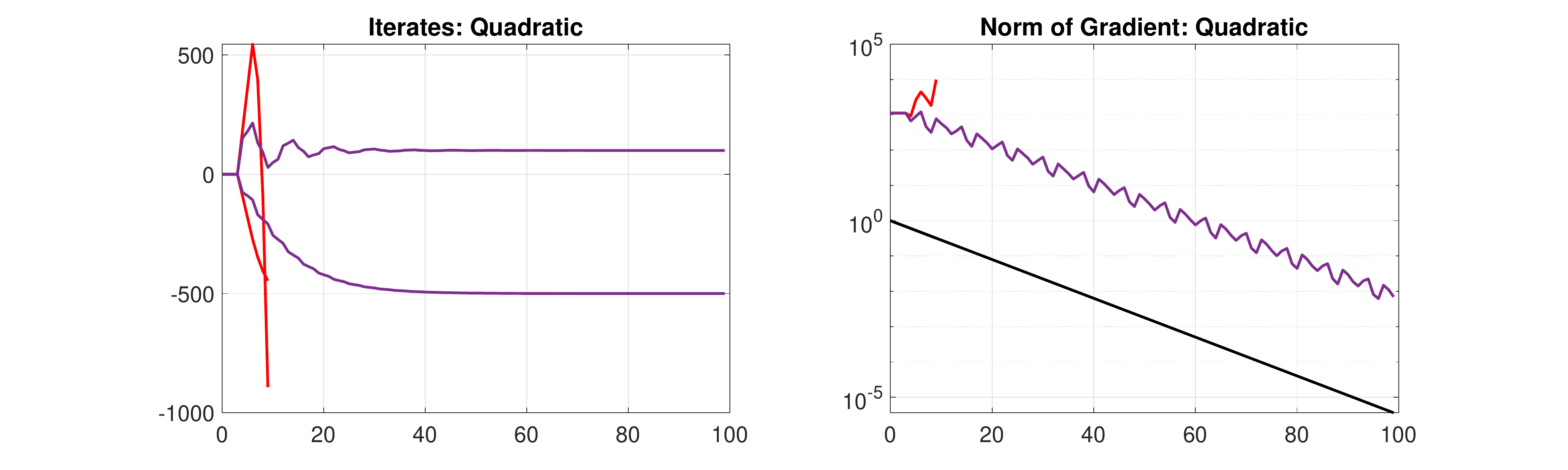}
\vspace*{-3ex}
\caption{Iterations for Case 4 with $m=1$, $L=10$ and a quadratic function
for optimal algorithm (purple) and gradient decent (red). The black line indicates
the optimal rate for synthesis.
}\label{figSim2a}
\end{figure}

\begin{figure}\center
\includegraphics[width=.48\textwidth,trim=120 15 100 0,clip=true]{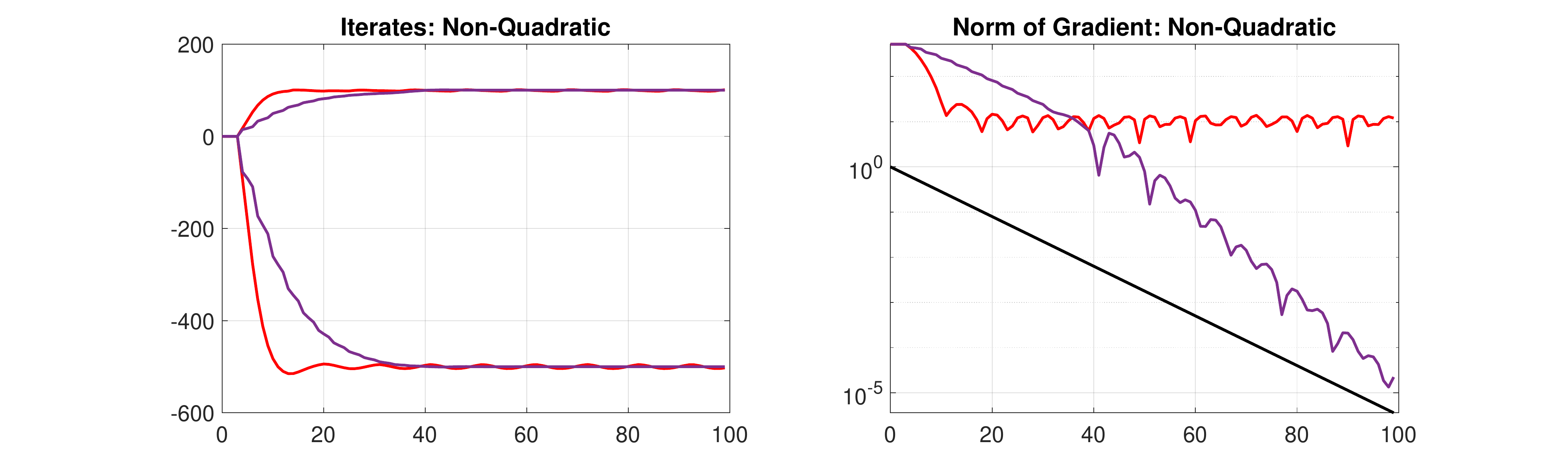}
\vspace*{-3ex}
\caption{
Plots corresponding to Fig.~\ref{figSim2a} for a non-quadratic function.}\label{figSim2b}
\end{figure}

Figs.~\ref{figSim2a} and \ref{figSim2b} depict the
iterates with the optimal algorithm in Case 4,
both for the quadratic and non-quadratic functions
$f(x)=h(x-b)$ with
$h(x)=Lx_1^2+mx_2^2$,
$h(x)=\frac{1}{2}(x_1^2+_2^2)+9\log(\exp(-x_1)+\exp(\frac{1}{3}x_1+x_2)+\exp(\frac{1}{3}x_1-x_2))$ \cite{Les22}
and $b=(100,-500)$ for
$m=1$ and $L=10$, respectively.
The optimal rates (depicted by the black line) are matched
in the quadratic case and give an upper bound for the non-quadratic function.
Gradient decent fails to converge in both cases.

\begin{figure}\center
\includegraphics[width=.48\textwidth,trim=120 15 100 0,clip=true]{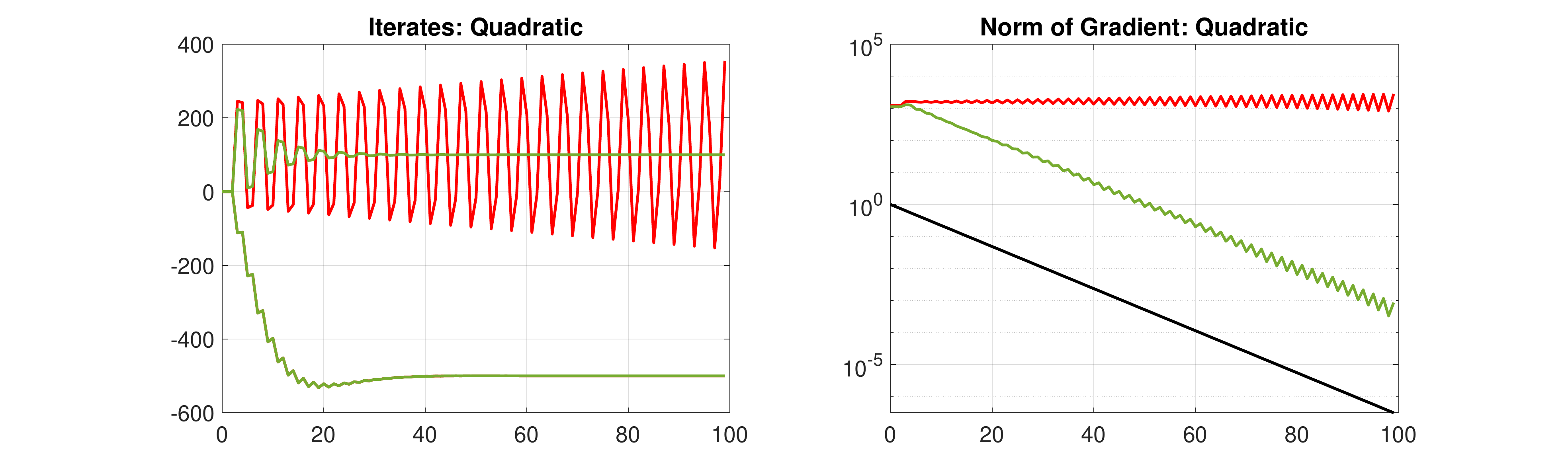}
\vspace*{-3ex}
\caption{
Iterations for Case 5 with optimal algorithm designed for $m=1$ and $L=10$
and employed for a quadratic function with $m=1$, $L=10$ (green)
as well as for $m=1$, $L=11$ (red), respectively. The black line indicates
the optimal rate for synthesis.}\label{figSim3a}
\label{figSim3}
\end{figure}


Finally, Fig.~\ref{figSim3a} reveals that optimal algorithms are
working well for the classes of functions they are designed for.
As expected from robust control, however, they can be sensitive to deviations from the assumptions,
as they lead to instability for the class $\c{S}_{1,11}$ that is only slightly larger than
$\c{S}_{1,10}$. All presented results can be reproduced with the software at \cite{code}.

\section{Conclusions}\label{Scon}
In this paper we have presented the full pipeline to design optimal optimization algorithms
or extremum controllers based on causal dynamic stability multipliers. A novel
parametrization of these filters overcomes technical assumptions as required in previous work.
Future work is devoted to incorporating inexact gradient information and to handling
anti-causal multipliers as well as performance
objectives in synthesis.





\begin{thebibliography}{10}
\providecommand{\url}[1]{#1}
\csname url@samestyle\endcsname
\providecommand{\newblock}{\relax}
\providecommand{\bibinfo}[2]{#2}
\providecommand{\BIBentrySTDinterwordspacing}{\spaceskip=0pt\relax}
\providecommand{\BIBentryALTinterwordstretchfactor}{4}
\providecommand{\BIBentryALTinterwordspacing}{\spaceskip=\fontdimen2\font plus
\BIBentryALTinterwordstretchfactor\fontdimen3\font minus
  \fontdimen4\font\relax}
\providecommand{\BIBforeignlanguage}[2]{{%
\expandafter\ifx\csname l@#1\endcsname\relax
\typeout{** WARNING: IEEEtran.bst: No hyphenation pattern has been}%
\typeout{** loaded for the language `#1'. Using the pattern for}%
\typeout{** the default language instead.}%
\else
\language=\csname l@#1\endcsname
\fi
#2}}
\providecommand{\BIBdecl}{\relax}
\BIBdecl

\bibitem{Nes18}
Y.~Nesterov, \emph{Lectures on Convex Optimization}, ser. Springer Optimization
  and Its Applications.\hskip 1em plus 0.5em minus 0.4em\relax Springer
  International Publishing, 2018, vol. 137.

\bibitem{Pol87}
B.~Polyak, \emph{Introduction to Optimization}.\hskip 1em plus 0.5em minus
  0.4em\relax Optimization Software, Inc., New York, 1987.

\bibitem{WanEli11}
J.~Wang and N.~Elia, ``A control perspective for centralized and distributed
  convex optimization,'' in \emph{Proceedings of the IEEE Conference on Decision and
  Control and European Control Conference, Orlando, FL, 2011}.

\bibitem{DueEbe12}
H.-B. D\"urr and C.~Ebenbauer, ``On a class of smooth optimization algorithms
  with applications in control,'' {IFAC} Proc. vol.~45, no.~17.\hskip 1em plus 0.5em minus
  0.4em\relax Elsevier {BV}, 2012, pp. 291--298.

\bibitem{LesRec16}
L.~Lessard, B.~Recht, and A.~Packard, ``{A}nalysis and {D}esign of
  {O}ptimization {A}lgorithms via {I}ntegral {Q}uadratic {C}onstraints,''
  \emph{SIAM Journal on Optimization}, vol.~26, no.~1, pp. 57--95, 2016.


\bibitem{Nes83}
Y.~Nesterov, ``A method for unconstrained convex minimization problem with the
  rate ofconvergence ${O}(\frac{1}{k^2})$.'' \emph{Doklady AN SSSR}, vol. 269,
  pp. 543--547, 1983, (In Russian; translated as SovietMath. Docl.).

\bibitem{ScoFre18}
B.~V. Scoy, R.~A. Freeman, and K.~M. Lynch, ``The fastest known globally
  convergent first-order method for minimizing strongly convex functions,''
  \emph{IEEE Control Systems Letters}, vol.~2, no.~1, pp. 49--54, 2018.



\bibitem{ZhoDoy96}
K.~Zhou, J.~Doyle, K. Glover, \emph{Robust and Optimal Control},
.\hskip 1em plus 0.5em minus
  0.4em\relax Prentice Hall, 1996.

\bibitem{SchWei11}
C.~Scherer, S.~Weiland,
``Linear matrix inequalities in control,'' in
  \emph{The Control Systems Handbook, Second Edition: Control System Advanced Methods},
  CRC Press, Chapter 24, pp. 1-30, 2011.


\bibitem{LesSei20}
L.~Lessard and P.~Seiler, ``Direct synthesis of iterative algorithms with
  bounds on achievable worst-case convergence rate,'' in \emph{2020 American
  Control Conference}, 2020, pp. 119--125.

\bibitem{MicSch21}
S.~Michalowsky, C.~Scherer, and C.~Ebenbauer, ``Robust and structure exploiting
  optimisation algorithms: an integral quadratic constraint approach,''
  \emph{International Journal of Control}, vol.~94, pp. 1--24, 2021.

\bibitem{GraEbe22}
D.~Gramlich, C.~Ebenbauer, and C.~W. Scherer, ``Synthesis of accelerated
  gradient algorithms for optimization and saddle point problems using
  {L}yapunov functions and {LMIs},'' \emph{Systems {\&} Control Letters}, vol.
  165, 2022.

\bibitem{TayHen16}
A.~B. Taylor, J.~M. Hendrickx, and F.~Glineur, ``Smooth strongly convex
  interpolation and exact worst-case performance of first-order methods,''
  \emph{Mathematical Programming}, vol. 161, no. 1-2, pp. 307--345, 2016.

\bibitem{TayDro22}
A.~B. Taylor and Y.~Drori, ``An optimal gradient method for smooth strongly convex
  minimization,'' \emph{Mathematical Programming}, vol. 199, no. 1-2, pp.
  557--594, 2022.

\bibitem{SchEbe21}
C.~Scherer and C.~Ebenbauer, ``Convex synthesis of accelerated gradient
  algorithms,'' \emph{{SIAM} Journal on Control and Optimization}, vol.~59,
  no.~6, pp. 4615--4645, 2021.

\bibitem{HolSch21c}
T.~Holicki and C.~W. Scherer, ``Algorithm design and extremum control: Convex
  synthesis due to plant multiplier commutation,'' in \emph{60th {IEEE}
  Conference on Decision and Control}, 2021, pp. 3249--3252.

\bibitem{Sch23}
C.~W. Scherer, ``Robust exponential stability and invariance guarantees with
  general dynamic {O}'{S}hea-{Z}ames-{F}alb multipliers,'' in \emph{Proc. IFAC
  World Congress, to appear}, 2023.

\bibitem{Bec17}
A.~Beck, \emph{First-Order Methods in Optimization}.\hskip 1em plus 0.5em minus
  0.4em\relax Society for Industrial and Applied Mathematics, 2017.

\bibitem{DesVid75}
C.~Desoer and M.~Vidyasagar, \emph{{F}eedback Systems: {I}nput-Output
  Approach}.\hskip 1em plus 0.5em minus 0.4em\relax London: Academic Press,
  1975.

\bibitem{MegRan97}
A.~Megretski and A.~Rantzer, ``{S}ystem analysis via {I}ntegral {Q}uadratic
  {C}onstraints,'' \emph{IEEE T. Automat. Contr.}, vol.~42, pp. 819--830, 1997.

\bibitem{WilBro68}
J.~Willems and R.~Brockett, ``Some new rearrangement inequalities having
  application in stability analysis,'' \emph{IEEE T. Automat. Contr.}, vol.~13,
  no.~5, pp. 539--549, 1968.

\bibitem{Zamfal68}
G.~Zames and P.~L. Falb, ``{S}tability conditions for systems with monotone and
  slope-restricted nonlinearities,'' \emph{SIAM Journal of Control}, vol.~6,
  pp. 89--109, 1968.

\bibitem{BocLes15}
R.~Boczar, L.~Lessard, and B.~Recht, ``Exponential convergence bounds using
  integral quadratic constraints,'' in \emph{54th {IEEE} Conference on Decision
  and Control ({CDC})}, 2015, pp. 7516--7521.

\bibitem{HuSei16}
B.~Hu and P.~Seiler, ``Exponential decay rate conditions for uncertain linear
  systems using integral quadratic constraints,'' \emph{{IEEE} Transactions on
  Automatic Control}, vol.~61, no.~11, pp. 3631--3637, 2016.

\bibitem{Sch22}
C.~W. Scherer, ``Dissipativity and integral quadratic constraints: Tailored
  computational robustness tests for complex interconnections,'' \emph{{IEEE}
  Control Systems Magazine}, vol.~42, no.~3, pp. 115--139, 2022.

\bibitem{GahApk94}
P.~Gahinet and P.~Apkarian, ``{A} linear matrix inequality approach to
  {$H_\infty$} control,'' \emph{Internat. J. Robust Nonlinear Control}, vol.~4,
  pp. 421--448, 1994.

\bibitem{IwaSke94}
T.~Iwasaki and R.~Skelton, ``{All} controllers for the general {${\cal
  H}_\infty$} control problem: {LMI} existence conditions and state space
  formulas,'' \emph{Automatica}, vol.~30, pp. 1307--1317, 1994.

\bibitem{Gah94}
P.~Gahinet, ``{A} new parametrization of {$H_\infty$} suboptimal controllers,''
  pp. 1031--1051, 1994.

\bibitem{Les22}
L.~Lessard, ``The analysis of optimization algorithms: A dissipativity
  approach,'' \emph{{IEEE} Control Systems Magazine}, vol.~42, no.~3, pp. 58--72, 2022.

\bibitem{AstMur09}
K.~J. Astr{\"o}m and R.~M. Murray, \emph{Feedback Systems: An Introduction for
  Scientists and Engineers}.\hskip 1em plus 0.5em minus 0.4em\relax Princeton
  Universety Press, 2009.

\bibitem{SouBha81}
E.~de~Souza and S.~Bhattacharyya, ``Controllability, observability and the
  solution of ${AX} - {XB} = {C}$,'' \emph{Linear Algebra and its Applications},
  vol.~39, pp. 167--188, 1981.

\bibitem{code}
A Matlab implementation of the results of this paper can be found at 
{\tt https://zenodo.org/badge/latestdoi/691960972}

\ls{
\bibitem{SchEbe23a}
C.W.~Scherer, C.~Ebenbauer and T.~Holicki, ``Optimization Algorithm Synthesis based on Integral Quadratic Constraints: A Tutorial,'' {\tt https://doi.org/10.48550/arXiv.2306.00565}, 2023
}{}

\end{thebibliography}

\ls{}{
\section{Relations to Classical Passivity Theory}\label{A0}

\begin{figure}\center
\scalebox{0.9}{
\begin{tikzpicture}[xscale=1,yscale=1]
\def\dl{4*\dn}
\node[sy2
] (f)  at (0,0) {$\nabla f_m\circ(\nabla f^L)^{-1}$};
\node[sy3,right= 1.5*\dl of f] (iir) {$\phi^{-1}$};
\node[sy2, below = .6*\dl of f] (d)   {$G_{m,L}$};
\node[sy3, below = .83*\dl of iir] (ir) {$\phi$};

\draw[<-] (f)--node[]{} (iir);
\draw[->] (f.west)-- ([xshift=-1.5*\dl] f.west) |- node[pos=.75]{$q$} (d.west) ;
\draw[<-] (iir.east)-- ([xshift=1.5*\dl] iir.east) |- node[swap,pos=0.75]{$r$}(ir.east);

\draw[->] (d)--node[]{$p$} (ir);
\end{tikzpicture}}
\vspace*{-1ex}
\caption{Classical loop with stability multiplier.
}\label{fig5}
\end{figure}
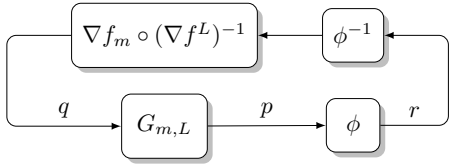

To explain the relation of Theorems~\ref{Trs} and \ref{Trs2} to the classical passivity theorem,
we consider the case $\rho=1$
and assume that $f^L$ is strongly convex. Then $\nabla f^L:\R^d\to\R^d$ has a global inverse $(\nabla f^L)^{-1}:\R^d\to\R^d$. Let us abbreviate the system \r{sysrhof} with \r{sya} and the filter \r{fild} as
$$
p=G_{m,L}(q)\te{and}r=\phi(p).
$$

If filtering the signals of the loop \r{nl}-\r{sy}
as in \r{psi}, a simple calculation leads to the
interconnection
\eql{looptrafo}{
q=\nabla f_m\circ(\nabla f^L)^{-1}(p)\te{and}p=G_{m,L}(q)
}
as depicted in Fig.~\ref{fig5} in case that $\phi$ is the static gain $1$.
Since the state-trajectories of the respective linear systems are identical,
stability of the original loop \r{nl}-\r{sy} is equivalent to that
of the transformed loop \r{looptrafo}.

Next note that \r{iqc0}-\r{psi} translate into the fact that the new map $\nabla f_m\circ(\nabla f^L)^{-1}$ is passive.
Moreover, \r{lmia} just expresses
that $-G_{m,L}$ is strictly passive.
Therefore, the passivity theorem
guarantees that the new loop \r{looptrafo} is asymptotically stable, and
Theorem~\ref{Trs} just boils down to this classical result applied to the transformed loop
\r{looptrafo}. Note that $\nabla f\mapsto \nabla f_m\circ (\nabla f^L)^{-1}$ can be
interpreted as taking maps in the sector $[m,L]$ into
passive ones. Moreover, strict passivity of $-G_{m,L}$ translates, by the KYP-Lemma, into the
classical circle criterion expressed in terms of the transfer function of the linear system \r{sy}.

If $\phi$ is nontrivial and dynamic, Lemma~\ref{Lpf} implies that $\phi$ and $\phi^{-1}$ are both stable. Hence, guaranteeing stability of the loop \r{nl}-\r{sy}
is equivalent to guaranteeing stability of the one in Fig.~\ref{fig5} involving the multiplier $\phi$ in a classical sense. Lemma~\ref{Lzfd} just means that
$q=\nabla f_m\circ(\nabla f^L)^{-1}(\phi^{-1}r)$
is passive for {\em any} multiplier $\phi$ subject to \r{lac}. Moreover,
Theorem~\ref{Trs2} based on \r{sycf} expresses that $-\phi G_{m,L}$ is strictly passive
for {\em some} multiplier $\phi$ subject to \r{lac}.
Hence, Theorem~\ref{Trs2} is an incarnation of the passivity theorem with causal multipliers
as addressed in detail in \cite[Chapter 9]{DesVid75}, with the additional feature of merging
the reduction of conservatism over the passivity theorem with a corresponding computational search over a convex family of multipliers.

All this illustrates the key ideas underlying the more powerful general dissipativity theory involving integral quadratic constraints, as exposed for example in \cite{Sch22}.

\section{Appendix: Proofs and an Auxiliary Result}

\subsection{Proof of Lemma~\ref{Ldi}}\label{AL2}
The proof starts with standard arguments in convex analysis.
Set $\al:=L-m>0$. Fix any $y\in\R^d$ and define
$$
g(x):=f_m(x)-\nabla f_m(y)^\T x\te{for}x\in\R^n.
$$
Clearly, $g$ is convex since it is an affine perturbation of the convex function $f_m$.
Due to
$\al q(x)-g(x)= Lq(x)-f(x)+\nabla f_m(y)^\T x,$ we conclude
that  $\al q-g$ is also convex. An application of the subgradient inequality for $\al q-g$ leads to
\eql{sg}{g(x+h)\leq g(x)+\nabla g(x)^\T h+\al q(h)}
for all $x,h\in\R^d$. Since $\nabla g(y)=\nabla f_m(y)-\nabla f_m(y)=0$, we conclude
$g(y)\leq g(x+h)$ and hence, with \r{sg}, that
$$
g(y)\leq g(x)+\nabla g(x)^\T h+\al q(h)\te{for all}x,h\in\R^d.
$$
The minimum of the convex quadratic function in $h$ on the right is easily calculated as
$g(x)-\frac{1}{\al}q(\nabla g(x))$.  This implies
$$
\al g(y)\leq \al g(x)-q(\nabla g(x))\te{for all}x\in\R^d.
$$
Simple rearrangements lead to \r{di}. Indeed,
with the definition of $g$ and $q(u-v)=q(u)+v^\T (v-u)-q(v)$ for $u=\nabla f_m(x)$, $v=\nabla f_m(y)$, we get
$$
\al f_m(y)-\al f_m(x)-\al v^\T(y-x)\leq -q(u)-v^\T(v-u)+q(v)
$$
and this gives
$$[\al f_m(y)-q(v)]-[\al f_m(x)-q(u)] \leq v^\T[\al y-v-(\al x-u)].$$
Since $\al y-v=(L-m)y-\nabla f(y)+my=Ly-\nabla f(y)$ and
$\al x-u=Lx-\nabla f(x)$, we finally get
\mun{[\al f_m(y)-q(\nabla f_m(y))]-[\al f_m(x)-q(\nabla f_m(x))]\leq\\ \leq
\nabla f_m(y)^\T[Ly-\nabla f(y)-(Lx-\nabla f(x))].}
This is \r{di}.
In case of $f\in\c{S}_{m,L}^0$ we infer $V(0)=0$ by its very definition.
Since $\nabla f_m(0)=0$, we can apply \r{di} for $u=0$ to infer that $V$ is indeed globally nonnegative.

\subsection{Proof of Theorem~\ref{Trs}}\label{ATrs}

Due to \r{lmia} and using standard rules for the Kronecker product, there exist some $\eps>0$ with
\mul{\label{disd}
[\bul]^\T\!\!\mat{cc}{\Xcl\ot \Id&0\\0&-(\Xcl\ot \Id)}\mat{cc}{\Acl\ot \Id&\Bcl\ot \Id\\I&0}+\\
+[\bul]^\T\!\!\mat{cc}{0&\Id\\\Id&0}\mat{cc}{\Ccla\ot \Id&\Dcla \ot \Id\\0&\Id}+
\mat{cc}{\eps I&0\\0&0}\cle 0.}
Now pick any trajectory of \r{nl}-\r{sy}. This can be transformed into one of \r{sysrho}-\r{nlrho}
and filtered by \r{psi} to generate a trajectory of \r{sysrhof}. Right-multiplying \r{disd} with $\col(\bar x_t,\bar q_t)$
and left-multiplying the transposed signal then leads to the dissipation inequality
$$
\bar x_{t+1}^\T(\Xcl\ots \Id)\bar x_{t+1}-\bar x_t^\T(\Xcl\ots \Id)\bar x_t+
2\bar q_t^\T\bar p_t+\eps\|\bar x_t\|^2\leq 0
$$
for all $t\in\N_0$. Now we exploit that the loop trajectories are as well related as in
\r{nlrho}. From \r{psi} and \r{gra} we infer
$\bar p_t=\bar F^L(t,\bar z_t)$, $\bar q_t=\bar F_m(t,\bar z_t)$, which shows
$\sum_{t=0}^{T-1}\bar q_t^\T\bar p_t\geq 0$ for all $T\in\N$  due to \r{iqc0}.
Summation of the dissipation inequality for $t=0,\ldots,T-1$ hence implies
\eql{dis2}{
\bar x_T^\T (\Xcl\ot \Id) \bar x_T-\bar x_0^\T (\Xcl\ot \Id)\bar x_0+
\eps\sum_{t=0}^{T-1} \|\bar x_t\|^2\leq 0
}
for all $T\in\N$. If $0<k_-<k_+$ denote the smallest and largest eigenvalues of $\c{X}\ot \Id$,
we conclude
\eql{dis3}{
k_-\|\bar x_T\|^2+\eps\sum_{t=0}^{T-1} \|\bar x_t\|^2\leq k_+\|\bar x_0\|^2\te{for all}T\in\N.
}
Since $\bar x_T=\rho^{-T}x_T$ and $\bar x_0=x_0$, we obtain \r{expsta} with $K=k_+/k_-$.
For $\rho=1$, we can conclude from \r{dis3} that $\eps\sum_{t=0}^\infty \|x_t\|^2<\infty$ and, hence,
$\lim_{t\to\infty}x_t=0$.

\subsection{Proof of Lemma~\ref{Lzf}}\label{AL4}
Fix $\nu\in\N$ and define $z_t:=\rho^t\bar z_t$ for $t\in\N_0$.
From \r{di} and \r{dis} we infer for all $t\in\N_0$ that
$$
V(z_t)-V(z_{t-\nu})\leq
\nabla f_m(z_t)^\T[\nabla f^L(z_t)-\nabla f^L(z_{t-\nu})],
$$
$$
V(z_t) \leq \nabla f_m(z_t)^\T\nabla f^L(z_t).
$$
The conic combination with the coefficients $\rho^{-2(t-\nu)}\geq 0$ and
$\rho^{-2t}-\rho^{-2(t-\nu)}=\rho^{-2t}(1-\rho^{2\nu})\geq 0$ gives
\mul{\label{h0}
\rho^{-2t}V(z_t)-\rho^{-2(t-\nu)}V(z_{t-\nu})\leq\\\leq
\nabla f_m(z_t)^\T[\rho^{-2t}\nabla f^L(z_t)-\rho^{-2(t-\nu)}\nabla f^L(z_{t-\nu})].
}
The right hand side equals
\mul{\label{hL1}
\rho^{-t}\nabla f_m(z_t)^\T[\rho^{-t}\nabla f^L(z_t)-\rho^\nu\rho^{-(t-\nu)}\nabla f^L(z_{t-\nu})]=\\=
\bar F_m(t,\bar z_t)^\T(\bar F^L(t,\bar z_t)-\rho^{\nu}\bar F^L(t-\nu,\bar z_{t-\nu})).
}
Now note for the left-hand side of \r{h0} that
\mun{
\sum_{t=0}^{T-1}\left(\rho^{-2t}V(z_t)-\rho^{-2(t-\nu)}V(z_{t-\nu})\right)=\\=
\sum_{t=0}^{T-1}\rho^{-2t}V(z_t)\ -
\!\!\!\sum_{t=-\nu}^{T-1-\nu}\!\!\rho^{-2t}V(z_{t})=
\!\!\!\sum_{t=T-\nu}^{T-1}\!\!\rho^{-2t}V(z_{t})\geq 0,
}
where we exploited $V(z_t)=0$ for $t<0$. Summation of \r{h0} and using \r{hL1} hence proves \r{iqc1}.

\subsection{Proof of Lemma~\ref{Lzfd}}\label{AL5}

Let $\mu_\nu:=-\la_\nu\rho^{-\nu}\geq 0$ for $\nu\in\N$ and
$\mu_0:=\sum_{\nu=1}^{\infty}\mu_\nu$.  Then
$\la_0-\mu_0=\la_0+\sum_{\nu=1}^\infty \la_\nu\rho^{-\nu}>0$.
Let us now conically combine \r{iqc0} and \r{iqc1} to infer
for any $T\in\N$ that
$$
0\leq
(\la_0-\mu_0)
\sum_{t=0}^{T-1}\bar q_t^\T\bar p_t+
\sum_{\nu=1}^\infty\mu_\nu
\sum_{t=0}^{T-1}\bar q_t^\T[\bar p_t-\rho^{\nu}\bar p_{t-\nu}].
$$
Since $\mu_0=\sum_{\nu=1}^{\infty}\mu_\nu$, this simplifies to
$$
0\leq
\la_0
\sum_{t=0}^{T-1}\bar q_t^\T\bar p_t
+\sum_{t=0}^{T-1}\bar q_t^\T\left[\sum_{\nu=1}^t(-\mu_\nu\rho^{\nu})\bar p_{t-\nu}\right],
$$
where we also exploited that
$\bar p_{t-\nu}=\bar F^L(t-\nu,\bar z_{t-\nu})=0$ for $\nu>t$.
Recalling the definition of $\mu_\nu$ leads to
$$
0\leq
\sum_{t=0}^{T-1}\bar q_t^\T
\left[\la_0\bar p_t+\sum_{\nu=1}^t\la_\nu\bar p_{t-\nu}\right]=
\sum_{t=0}^{T-1}\bar q_t^\T\bar r_t
$$
which was to be shown.
\subsection{Proof of Lemma~\ref{Lpf}}\label{AL6}

We start by proving the equivalence of the
first properties in \r{iirc} and \r{lac} (for \r{fil}).
Indeed, \r{iirc} implies $\la_\nu\leq 0$ for $\nu=0,\ldots,l-1$,
which gives $\Cfla\c{K}(\Af,\Bf)\leq 0$. Conversely,
$\Cfla\c{K}(\Af,\Bf)\leq 0$ implies $\Cfla\Af^\nu\Bf\leq 0$ for $\nu=0,\ldots,l-1$.
By the Cayley-Hamilton theorem, we note that
\eql{ch}{
\Cfla\Af^{l+\mu}\Bf=\sum_{j=0}^{l-1}(-\al_j)\Cfla\Af^{j+\mu}\Bf\te{for all}\mu\in\N_0.
}
Since
$\al_j\leq 0$ for $j=0,\ldots,l-1$, we conclude from \r{ch}
for $\mu=0$ that  $\Cfla\Af^l\Bf\leq 0$.
An induction step based on \r{ch} for $\mu\in\N$ then proves
$\Cfla\Af^{\nu}\Bf\leq 0$ for all $\nu=l+1,l+2,\ldots.$

Since $\Af$ has all its eigenvalues in $\D_\rho$, we note next that
\eql{hfil}{
\sum_{\nu=0}^\infty (\Cfla\Af^\nu\Bf) \z^{-(\nu+1)}=\Cfla(\z I-\Af)^{-1}\Bf
}
if $|\z|\geq\rho$.
Using \r{hfil} for $\z=\rho$ shows that the second conditions in \r{iirc}  and \r{lac} (for \r{fil}) are as well equivalent.

Finally, assume that \r{lac} holds.
Based on \r{hfil} for $\z=\rho>0$, we can as well infer $\Cfla(\rho I-\Af)^{-1}\Bf\leq 0$ and thus $\Dfla>0$.
Then let $\la \in\C$ be an eigenvalue of $\Af-\Bf\Dfla^{-1}\Cfla$ with $|\la |\geq\rho$.
Since $\al(\la )\neq 0$, we infer that $\Dfla+\Cfla(\la  I-\Af)^{-1}\Bf=0$. Again with \r{hfil}, we infer
\mun{
0=|\Dfla+\Cfla(\la  I-\Af)^{-1}\Bf|\geq \\\geq
\Dfla-\sum_{\nu=1}^\infty |\Cfla\Af^\nu\Bf||\la^{-\nu}|\geq
\Dfla-\sum_{\nu=1}^\infty |\Cfla\Af^\nu\Bf|\rho^{-\nu}= \\=
\Dfla+\sum_{\nu=1}^\infty (\Cfla\Af^\nu\Bf)\rho^{-\nu}=\Dfla+\Cfla(\rho I-\Af)^{-1}\Bf>0.
}
This contradiction concludes the proof.

\subsection{Proof of Theorem~\ref{Trs2}}\label{ATrs2}

We follow the proof of Theorem~\ref{Trs}. This leads to a trajectory of \r{sysrhof},
but now for the matrices \r{sycf} and with the state-trajectory $\col(\xi_t,\bar x_t)$
comprising both the one of the filter \r{fild} and the system \r{sysrho}.
Since $\xi_0=0$, we conclude as earlier that
$$
\mat{c}{\xi_t\\\bar x_{t}}^\T\!\!
\underbrace{\mat{cc}{\Xclf&\Xclfs\\\Xclsf&\Xcls}}_{\Xcl\ot\Id}
\mat{c}{\xi_t\\\bar x_{t}}
+\eps\left\|\mat{c}{\xi_t\\\bar x_t}\right\|^2
\leq \bar x_{0}^\T\Xcls\bar x_{0}
$$
holds for all $t\in\N_0$, where
$\Xcl\ot\Id$ is partitioned according to $\Acl\ot\Id$ in \r{sycf}.
It remains to observe that the left-hand side can be bounded from below
by $\bar x_t^\T(\Xcls-\Xclsf\Xclf^{-1}\Xclfs)\bar x_t+\eps\|\bar x_t\|^2$ and that the appearing Schur-complement is positive definite. This permits to conclude the proof as earlier.

\subsection{Proof of Theorem~\ref{Tsyn}}\label{ATsyn}

Proof of $(a)\Rightarrow (b)$.
By a slight perturbation of $(\Cfla,\Dfla)$, we can make sure that
\eql{per}{
\arr{c}{
\text{$(\Af,\Cfla)$ is observable and no eigenvalue of $\At$ is}\\
\text{a zero of $\be(\z):=[\Dfla+\Cfla(\z I-\Af)^{-1}\Bf]\al(\z)$,}
}}
while  \r{lmia}, \r{lac} still hold true. Then all assumptions in
Lemma~\ref{Lcom} are satisfied for \r{subsys}. Note that
$\be(\z)=\Dfla\al(\z)+\Cfla\col(1,\z,\ldots\z^{l-1})$. Hence, equations
1)-4) in Lemma~\ref{Lcom} match with those formulated in the theorem,
and \r{Tra}
is invertible. By \r{lac} and Lemma~\ref{Lpf}, we conclude that
\eql{fii}{
\Af^i:=\Af-\Bf\Dfla^{-1}\Cfla\te{is Schur.}
}
Moreover, \r{tracom} in Lemma~\ref{Lcom} implies
\eql{h1}{
\mat{c|cc}{
\Tla\Ah \Tla^{-1} &\Tla\Bh_1  &\Tla\Bh\hl
\Chla\Tla^{-1}    &\Dhla       &\Ehla}=
\mat{cc|ccc}{
\Af        &0       &\Tf\Bh_1   &\Bf \\
\Bt\Cfla   &\At     &\Ttla\Bh_1 &\Bt\Dfla\hl
\Et\Cfla   &\Ctz   &\Dhla      &\Ehla }.
}
In the sequel, we also exploit the obvious relation
\eql{h2}{
\mat{cc}{
\Af         \\
\Bt\Cfla   \hl
\Et\Cfla   }=
\mat{cc|ccc}{
\Af^i  \\
0      \hl
0      }+
\mat{ccccc}{
\Bf        \\
\Bt\Dfla  \hl
\Et\Dfla  }\Dfla^{-1}\Cfla.
}

Our first goal is to determine some annihilators $\Uh$ and $\Vh$ with \r{ann2}.
By inspection, $\Uh$ as given in the theorem already has this property. For convexification, some specific $\Vh$ needs to be constructed on the basis of $V$ in \r{ann}. To this end, we note that
\r{ann} and $\Dfla\neq 0$ imply
$$\arc
\ker\mat{cc}{(\Tla\Bh)^\T &\bm{\h{E}}^\T}\!=\!
\ker\mat{ccc}{\Bf^\T&\Dfla\Bt^\T &\Dfla\Et^\T}
\supset\im\mat{c}{0\\V_2^\T\hl V_3^\T}
$$
if partitioning the columns of $V$ as $\mat{cc}{V_2&V_3}$ accordingly.
Thus, we can determine $\Vf$, $\Va$, $\Vb$ such that
\eql{h3}{
\mat{cc}{\Vf^\T&0\\ \Va^\T&V_2^\T\\\Vb^\T &V_3^\T}
\te{is a basis of}
\ker\mat{cc}{\Bh^\T\Tla^\T&\bm{\h{E}}^\T}.
}
As a consequence, \r{ann2} is satisfied with the choice
\eql{ann3}{
\Vh:=
\mat{cc|c}{\Vf&\Va&\Vb \\ 0&V_2&V_3}
\mat{cc}{\Tla &0\\0&I}.
}
By construction, we note that the block $\Vf$ has full row rank.

By our preparatory remarks, it is then guaranteed that
that there exist symmetric matrices  $X$ and $\Yh$
with \r{lmip}-\r{lmisc}.

Let us now introduce
\eql{Yhelp}{
Y:=\Tla \Yh\Tla^\T=\mat{cc}{Y_{\rm f}&\bul \\\bul &\Yt}.
}
A congruence transformation of \r{lmisc} then leads to
\eql{lmich}{
\mat{ccc}{Y&\Tla\\\Tla^\T&X}=
\mat{ccc}{Y_{\rm f}&\bul&\Tf\\\bul&\Yt&\Ttla\\\Tf^\T&\Ttla^\T&X}
\cg 0.
}
Canceling the first block row and column shows \r{lmic}.

In moving towards \r{lmid}, we note that \r{lmisd} is equivalent to
\eql{lmidh}{
\underbrace{\Vh\mat{cc|cc}{-\Tla^{-1}&\Ah\Tla^{-1}&0&\Bh_1\\0& \Chla\Tla^{-1}&-1&\Dhla}}_F
\mat{cc|cc}{Y&0&0&0\\0&-Y&0&0\hl 0&0&0&1\\0&0&1&0}\bul^\T\cg 0.
}
Due to \r{ann3}, the product of the two matrices on the left equals
$$
F=\mat{c}{F_1\\F_2}=\mat{cc|c}{\Vf&\Va&\Vb \\ 0&V_2&V_3}
\mat{cc|cc}{-I&\Tla\Ah\Tla^{-1}&0&\Tla\Bh\\0& \Chla\Tla^{-1}&-1&\Dhla}.
$$
If we recall \r{h1} and exploit \r{h2} with \r{h3}, this reads as
\eql{Lhelp}{
\mat{c}{F_1\\F_2}=\mat{cc|c}{\Vf&\Va&\Vb \\ 0&V_2&V_3}
\mat{cc|cc|ccc}{
-I&0   &\Af^i   &0     &0 &\Tf\Bh_1    \\
0 &-I  &0       &\At   &0 &\Ttla\Bh_1 \hl
0 & 0  & 0      &\Ctz &-1&\Dhla }.
}
Then
$F_2=
V
\smat{cc|cc|ccc}{
0 &-I  &0       &\At   &0 &\Ttla\Bh_1 \hl
0 & 0  & 0      &\Ctz &-1&\Dhla }.
$
As a consequence, after canceling the first block row and column of
\r{lmidh} and if recalling \r{Yhelp}, we arrive at  \r{lmid}.

Proof of $(b)\Rightarrow (a)$. Suppose that \r{lac},
\r{lmip} and \r{lmid}-\r{lmic} are feasible.
By perturbing $(\Cfla,\Dfla)$ if necessary, we can again
assume that \r{per} holds true.
Due to \r{fii}, we can solve the Stein equation
$\Yf-\Af^i\Yf(\Af^i)^\T=I$ and note that
$\Yf\cg 0.$
Motivated by \r{Yhelp}, we define
\eql{Yhelp2}{\Yh:=\Tla^{-1}\mat{cc}{\ga\Yf&0 \\0 &\Yt}\Tla^{-\T}}
with some still to-be-determined scalar $\ga>0$. For this choice of $\Yh$ and in view of \r{Lhelp},
the inequalities \r{lmich}-\r{lmidh} read as
\eql{help0}{
\mat{ccc}{\ga\Yf&S_{12}\\S_{21}&S_{22}}\cg 0\te{and}
\mat{ccc}{\ga\Vf\Vf^\T+R_{11}&R_{12}\\R_{21}&R_{22}}\cg 0,
}
respectively, where the $R$- and $S$-blocks are independent from $\ga$.
Moreover, \r{lmic} and \r{lmid} just guarantee that $S_{22}\cg 0$ and $R_{22}\cg 0$.
Since $\Yf\cg 0$ and $\Vf\Vf^\T\cg 0$ (because $\Vf$ has full row rank),
we can choose $\ga$ sufficiently large in order to enforce the validity of
\r{help0}, and thus of \r{lmich}-\r{lmidh} with $Y=\diag(\ga\Yf,\Yt)$. Again,
these translate back into \r{lmisd}-\r{lmisc}. In summary, we have shown
the existence of $X$ and $\Yh$ to guarantee \r{lmip}-\r{lmisc}, which
in turn implies the existence of a controller as was to be shown.

\subsection{On Commuting Transfer Functions}\label{Scom}

\renewcommand{\Al}{A}
\renewcommand{\Bl}{B}
\renewcommand{\Cl}{C}
\renewcommand{\Dl}{D}

\renewcommand{\Ar}{\h A}
\renewcommand{\Br}{\h B}
\renewcommand{\Cr}{\h C}
\renewcommand{\Dr}{\h D}

\lemma{\label{Lcom}
For SISO transfer functions with realizations
$$
\Gl=\mas{c|c}{\Al&\Bl\hl \Cl &\Dl}\te{and}\Gr=\mas{c|c}{\Ar&\Br\hl \Cr &\Dr },
$$
let $\all$ and $\alr$ denote the characteristic polynomials of $\Al\in\R^{\nl\times \nl}$, $\Ar\in\R^{\nr\times\nr}$, respectively, and define the  polynomial $\bel:=\Gl\all$ of degree at most $\nl$.

Let $(A,B)$ and $(A,C)$ be controllable and observable, respectively, and
suppose that $(\Al,\Ar)$ have no common eigenvalues and
$(\bel,\alr)$ and $(\all,\Gr)$ no common zeros in $\C$.

Then the two standard realizations of the products in
$\Gl\Gr=\Gr\Gl$ are related by a state-coordinate change as
\eql{tracom}{
\mas{cc|ccc}{
\Al  &\Bl \Cr   &\Bl \Dr     \\
0    & \Ar     &\Br             \hl
\Cl  &\Dl  \Cr  &\Dl \Dr }
\stackrel{\arraycolsep.1ex\smat{c|c}{L^{-1}&-L^{-1}K\hl NL^{-1}&M\!-\!NL^{-1}K}}{\longrightarrow}
\mas{cc|ccc}{
\Al         &0              &\Bl   \\
\Br\Cl        &\Ar            &\Br\Dl  \hl
\Dr\Cl        &\Cr            &\Dr\Dl  }
}
with the solutions of the matrix equations
\enu{
\item $\Al K-K\Ar+\Bl\Cr=0$,
\item $L\c{K}(\Al,\Bl)=\c{K}(\Al,\Bl \Dr-K\Br),$
\item $M\all(\Ar)=\bel(\Ar),$
\item $\Ar N-N\Al+\Br\Cl=0.$
}
}

{\bf Proof.} {\em Preparation.}
Note that $\Gl=\Dl+\ga\all^{-1}$ for $\ga:=\bel-\Dl\all$.
With the coefficients of
$\all(\z)=\all_0+\cdots+\all_{\nl-1}\z^{\nl-1}+\z^\nl$ and
$\ga(\z)=:\ga_0+\ga_1\z+\cdots+\ga_{\nl-1}\z^{\nl-1}$
and next to the companion matrix $C_\all$ of $\all$,
let us introduce
$$\renewcommand{\arraystretch}{.9}
H_\all:=
\mat{ccccc}{
\all_1&\cdots&\all_{n-1}&1 \\[-1ex]
\vdots&&\vdots&\vdots\\
\all_{n-1}&\cdots&0&0\\
1&\cdots&0&0},\
B_\ga:=\mat{c}{\ga_0\\[-1ex]\vdots\\\ga_{\nl-2}\\\ga_{\nl-1}}.
$$
With the polynomial vectors $r_\nl(\z)=\mat{cccc}{\z^0&\z^1&\cdots&\z^{n-1}}$ and $c_\nl(\z)=r_\nl(\z)^T$,
the Kalman controllability and observability matrices of $(\Al,\Bl)$, $(\Al,\Cl)$ can be directly expressed as
$r_{\nl}(\Al)(I_\nl\otimes \Bl)$, $(I_\nl\otimes\Cl)c_{\nl}(\Al)$,
respectively. Then $\Sl:=r_{\nl}(\Al)(I_\nl\otimes \Bl)H_\all$ is invertible and guarantees
 \cite{AstMur09}
\eql{tracon}{
\Sl^{-1}\Al \Sl=C_\all,\ \
\Sl^{-1}\Bl=e_\nl\te{and}
\Cl\Sl=B_\ga^T.}
Analogously, $\Tl:=H_\all(I_\nl\ot \Cl)c_{\nl}(\Al)$ is invertible with
\eql{tracon2}{
\Tl\Al \Tl^{-1}=C_\all^T,\ \ \Tl\Bl=B_\ga\te{and}\Cl\Tl^{-1}=e_\nl^T.
}

{\em Step 1.} Since $\Al$ and $\Ar$ have no common eigenvalues,
the solutions $K$ and $N$ of the Sylvester equations 1) and 4) exist and are unique.
With these two equations we get
\eql{t1}{\arraycolsep.4ex
\mat{cc|ccc}{
\Al      &\Bl \Cr       &\Bl \Dr     \\
0        & \Ar          &\Br             \hl
\Cl      &\Dl  \Cr      &\Dl \Dr }
\!\!\stackrel{\arraycolsep.1ex\smat{cc}{I&-K\\0&I}}\longrightarrow\!\!
\mat{cc|ccc}{
\Al      &0               &\Bl \Dr \!-\!K\Br    \\
0        &\Ar             &\Br             \hl
\Cl      &\Dl\Cr\!+\!\Cl K    &\Dl \Dr },
}
\eql{t2}{\arraycolsep.5ex
\mat{cc|ccc}{
\Al        &0            &\Bl   \\
0          &\Ar          &\Br\Dl\!-\!N\Bl  \hl
\Dr\Cl\!+\!\Cr N &\Cr        &\Dr\Dl  }
\!\!\stackrel{\arraycolsep.1ex\smat{cc}{I&0\\N&I}}\longrightarrow\!\!
\mat{cc|ccc}{
\Al         &0            &\Bl   \\
\Br\Cl      &\Ar          &\Br\Dl  \hl
\Dr\Cl      &\Cr          &\Dr\Dl  }.}
Since $\Gl\Gr=\Gr\Gl$, the transfer functions defined by both sides of \r{t1} and \r{t2} are identical
and decompose as in
\mul{\label{id1}
\mas{c|ccc}{
\Al      &\Bl \Dr -K\Br  \hl
\Cl      &\Dl \Dr        }
+
\mas{c|ccc}{
\Ar             &\Br     \hl
\Dl\Cr+\Cl K    &\Dl \Dr }=\\=
\mas{c|ccc}{
\Al          &\Bl   \hl
\Dr\Cl+\Cr N &\Dr\Dl  }
+
\mas{c|ccc}{
\Ar          &\Br\Dl-N\Bl  \hl
\Cr        &\Dr\Dl  }.}

{\em Step 2.}
The goal is to prove that $L^{-1}$ is invertible and
\eql{rea1}{
\mat{c|ccc}{
\Al        &\Bl \Dr -K\Br    \hl
\Cl        &\Dl \Dr }
\stackrel{\arraycolsep.1ex L^{-1}}\longrightarrow
\mat{c|ccc}{
\Al            &\Bl   \hl
\Dr\Cl+\Cr N   &\Dr\Dl  }.
}
By \r{id1} and since $\Al$ and $\Ar$ share no eigenvalues,
the system's transfer functions on the left and on the right of \r{rea1} are indeed
identical. Let us now prove, respectively,  that
$$
(\Al,\Bl \Dr-K\Br),(\Al,\Dr\Cl+\Cr N)
\text{\ is controllable, observable.}
$$
Indeed, suppose that $x^*\Al=\la x^*$ and $x^*(\Bl \Dr-K\Br)=0$. Then 1) gives $x^*K(\la I-\Ar)+x^*\Bl \Cr=0$ and thus
$x^*K=-x^*\Bl \Cr(\la I-\Ar)^{-1}.$ This implies $0=x^*(\Bl \Dr-K\Br)=x^*\Bl(\Dr+\Cr(\la I-\Ar)^{-1}\Br)=x^*\Bl \Gr(\la)$.
Since $\all(\la)=0$, we infer $\Gr(\la)\neq 0$ and thus $x^*\Bl =0$. This shows $x^*=0$, just because $(\Al,\Bl)$ is controllable. Hence $(\Al,\Bl \Dr-K\Br)$ is controllable.
The Sylvester equation 4) for $N$ permits to similarly prove that
$(\Al,\Dr\Cl+\Cr N)$ is observable.

As a conclusion, the realizations of the identical transfer functions in
\r{rea1} are minimal. This guarantees the existence of $L$ with \r{rea1},
and $L^{-1}$ satisfies $L^{-1}\c{K}(\Al,\Bl \Dr-K\Br)=\c{K}(\Al,\Bl)$ \cite{AstMur09}.
In summary, we have shown that 2) has a unique invertible solution $L$ which guarantees \r{rea1}.

{\em Step 3.}
Since no eigenvalue of $\Ar$ is a zero of
$\all$ and $\bel$, both $\all(\Ar)$ and
$\bel(\Ar)$ are invertible. Therefore, $M$ as in 3) exists
and is unique and invertible.

Now we prove that
\eql{rea2}{
\mat{c|ccc}{
\Ar           &\Br             \hl
\Dl\Cr+\Cl K  &\Dl \Dr }
\stackrel{\arraycolsep.1ex M}\longrightarrow
 \mat{c|ccc}{
\Ar             &\Br\Dl-N\Bl  \hl
\Cr             &\Dr\Dl        }.
}
By its very definition, $M$ commutes with $\Ar$, which implies $M\Ar M^{-1}=\Ar$.
To prove $\Cr M=\Dl\Cr+\Cl K$, we observe
(see e.g. \cite[Proof of Lemma 4]{SouBha81}) that 1) implies
$$
\all(\Al)K-K\all(\Ar)+r_{\nl}(\Al)
(I_{\nl}\otimes \Bl)H_\all
(I_{\nl}\otimes \Cr)
c_\nl(\Ar)=0.
$$
If using $\all(\Al)=0$, we infer with $\Sl$ from \r{tracon} that
\eql{id2}{
K=
\Sl
(I_{\nl}\otimes \Cr)
c_\nl(\Ar)\all(\Ar)^{-1}.
}
Since $\Cr$ and $\Cl$ are rows, we have
$\Cl \Sl(I_{\nl}\otimes \Cr)=
(\Cl \Sl\otimes 1)(I_{\nl}\otimes \Cr)=(\Cl \Sl\otimes \Cr)=
(1\otimes\Cr)(\Cl \Sl\otimes I_n)=\Cr(\Cl \Sl\otimes I_n)$.
By using \r{tracon}, we hence conclude from \r{id2} that
$$
\Cl K=\Cr(B_\ga^T\ot I_n)c_\nl(\Ar)\all(\Ar)^{-1}=\Cr\ga(\Ar)\all(\Ar)^{-1}.
$$
Finally, recall $\bel=\Dl\all+\ga$ and thus
$\bel(\Ar)=\Dl\all(\Ar)+\ga(\Ar)$.
By 3), this leads to
$M=\Dl I_n+\ga(\Ar)\all(\Ar)^{-1}$
and then
$\Cr M= \Cr\Dl+\Cr\ga(\Ar)\all(\Ar)^{-1}=\Dl\Cr +\Cl K$,
as to be shown.

In complete analogy, 4) implies
with $\Tl$ from \r{tracon2} that
$
N=-\all(\Ar)^{-1}r_{\nl}(\Ar)(I_{\nl}\ot \Br)\Tl.
$
This leads to
$$
N\Bl=-\all(\Ar)^{-1}r_{\nl}(\Ar)(B_\ga\ot I_n)\Br=-\all(\Ar)^{-1}\ga(\Ar)\Br.$$
Then $M\Br=\Dl \Br+\all(\Ar)^{-1}\ga(\Ar)\Br=\Br\Dl-N\Bl$.

{\em Step 4.} Continuing \r{t1} with \r{rea2} and \r{rea1} leads to
\mun{
\arraycolsep.5ex
\mat{cc|ccc}{
\Al      &\Bl \Cr       &\Bl \Dr     \\
0        &\Ar           &\Br             \hl
\Cl      &\Dl  \Cr      &\Dl \Dr }
\stackrel{\arraycolsep.1ex\smat{cc}{I&-K\\0&M}}\longrightarrow
\mat{cc|ccc}{
\Al      &0               &\Bl \Dr -K\Br    \\
0        &\Ar             &\Br\Dl-N\Bl      \hl
\Cl      &\Cr             &\Dl \Dr }\\
\stackrel{\arraycolsep.1ex\smat{cc}{L^{-1}&-L^{-1}K\\0&M}}\longrightarrow
\mat{cc|ccc}{
\Al           &0     &\Bl          \\
0             &\Ar   &\Br\Dl-N\Bl  \hl
\Dr\Cl+\Cr N  &\Cr   &\Dl \Dr      }.
}
In combination with \r{t2}, the matrix
$$
\mat{cc}{I&0\\N&I}\mat{cc}{L^{-1}&-L^{-1}K\\0&M}=
\mat{ccc}{L^{-1}&-L^{-1}K\\NL^{-1}&M-NL^{-1}K}
$$
indeed guarantees \r{tracom}. This latter definition
also shows that the transformation matrix is invertible.\epro
}

\addtolength{\textheight}{-12cm}   

\end{document}